\documentclass[preprint,12pt,3p]{elsarticle}

\usepackage{amssymb}
\usepackage{amsmath,amsthm}
\usepackage{amsfonts}
\usepackage{enumerate}
\usepackage{graphicx}
\usepackage{comment}
\usepackage{enumitem}
\usepackage{cancel}
\usepackage{booktabs}
\usepackage{pgfplots}
\usepackage{subcaption}
\usepackage{float}
\usepackage{bm}
\pgfplotsset{compat=1.18}

 \usepackage{lineno}
 \usepackage[utf8]{inputenc}
\usepackage[figuresleft]{rotating}

\usepackage{dcolumn}
\usepackage{siunitx}
\usepackage{algorithm}
\usepackage{algpseudocode}
\usepackage{multirow}
\usepackage{xcolor}
\usepackage[normalem]{ulem}
\usepackage{graphicx}               
\usepackage{caption}                
\usepackage{subcaption}             
\usepackage{tikz}
\usepackage{comment}
\usetikzlibrary{graphs}
\usepackage{color}
\definecolor{orange}{rgb}{1,0.3,0.0}
\definecolor{violet}{rgb}{0.54, 0.17, 0.89}
\definecolor{aqua}{rgb}{0.0, 1.0, 1.0}

\newcommand{\ABC}{LRP-DP}

\newtheorem{rem}{Remark}
\newtheorem{proposition}{Proposition}

\journal{Transportation Research Part B}

\begin{document}

\begin{frontmatter}

\title{Handling Customer Preferences in Last-Mile Location-Routing: A Bilevel Optimization Approach}


\author[label1]{ Jos\'e Emmanuel Gómez-Rocha}
\address[label1]{Tecnologico de Monterrey, Escuela de Ingeniería y Ciencias, Mexico}
\ead{emmanuel.gr@tec.mx}

\author[label1]{ Jos\'e-Fernando Camacho-Vallejo \corref{cor1}}
\ead{fernando.camacho@tec.mx}
\cortext[cor1]{I am corresponding author}

\author[label2]{Elena Fern\'{a}ndez}
\address[label2]{Departamento de Estadística e Investigación Operativa, Universidad de C\'adiz, Spain}
\ead{elena.fernandez@uca.es}

\begin{abstract}
This paper introduces the Location-Routing Problem with Delivery Preferences,  
a last-mile delivery problem with service flexibility that integrates collection-point activation and routing decisions while accounting for customer choice. The company determines which collection points to activate, which customers to offer home delivery, and the corresponding delivery route. Customers subsequently decide whether to accept the proposed home-delivery service or collect their parcels from an activated collection point according to ranked preferences over delivery time intervals and collection points.\\
This hierarchical interaction is modeled through a bilevel optimization framework, which can be formulated as a mixed-integer bilevel program that incorporates compensation mechanisms for service inconvenience. Two equivalent single-level reformulations that exploit the structure of the customers' problems are developed. In addition, a two-phase heuristic is proposed to generate feasible solutions that can be used as warm starts for the formulations. The computational experiments run to compare the performance of the formulations, show that the warm-start strategy becomes particularly effective as the instance size increases. \\
The comparison with a benchmark model that ignores customer preferences indicates that explicitly accounting for customer choice alters the resulting location-routing configuration while improving customer preference satisfaction at a small increase in the company's overall cost. These findings highlight the importance of jointly considering operational efficiency, service flexibility, and customer autonomy in last-mile delivery planning.
\end{abstract}


\end{frontmatter}


\section{Introduction}\label{sec:Intro}
The rapid growth of e-commerce has profoundly reshaped urban freight distribution, placing significant pressure on last-mile logistics systems. Increasing order volumes, tighter delivery expectations, and greater consideration of customer opinions regarding service quality have made distribution planning considerably more challenging \citep{risberg2023systematic}. Failed home deliveries, traffic congestion, and rising operational costs further complicate last-mile logistics processes. In response, companies are progressively incorporating alternative delivery options, such as parcel collection points, to enhance flexibility and reduce costly redelivery attempts \citep{ma2022parcel, ozyavas2025designing}. Integrating these alternatives into routing decisions is therefore of significant practical and strategic importance.\\

At the same time, customers have become active participants in the delivery process. Unlike in the past, they now express explicit preferences regarding delivery time windows and pickup locations \citep{karli2024innovative}. The flexibility associated with these preferences can have important operational implications. For instance, \citet{yildiz2020pricing} show that additional delivery-time flexibility can generate substantial savings in distribution costs. As noted by \cite{cebeci2023changing}, ignoring customer behavior can lead to inefficient routing plans, service dissatisfaction, and reputational damage. Consequently, modern last-mile systems must be designed not only from an operational efficiency perspective but also with explicit consideration of customer-driven service constraints.\\

Motivated by these challenges, this paper introduces the Location-Routing Problem with Delivery Preferences (\ABC). The \ABC\ is 
a business-to-customer delivery problem in which an e-commerce company operates a central depot and may activate a set of collection points. The depot serves exclusively as the origin and destination of delivery routes, whereas collection points act as intermediate facilities where customers can retrieve their parcels as an alternative to home delivery. The service route visits all activated collection points and a subset of customers who are offered home delivery. Upon purchasing a product, each customer specifies an ideal time window for home delivery. Consequently, the company must design routes that account for these preferred service times and determine which customers will be offered home delivery.\\

Customers respond to the company's decisions according to their preferences as follows. Customers who are offered home delivery may accept home service even if the proposed service time differs from their ideal time. Alternatively, they may prefer to pick up their parcels at an activated collection point. On the other hand, customers who are not offered home delivery select their most preferred activated collection point for parcel pickup.\\

When home delivery does not take place at the customer's ideal service time, the company incurs compensation payments reflecting its estimate of the resulting service deterioration. These compensations are determined based on the customer information available to the company. For customers ultimately receiving home delivery, compensation depends on the deviation from their ideal service time. For customers served through a collection point, compensation is based on the distance to the nearest activated collection point, regardless of the collection point ultimately selected according to their preferences. This distinction separates the compensation attributable to the network designed by the company from any additional travel resulting from the customer's individual choice, which is independent of the service quality offered by the company. Therefore, customer reactions are incorporated into the company's routing and collection-point activation decisions to effectively design the last-mile logistics network.\\ 

The company aims to minimize total costs, including collection-point activation costs, routing costs, and customer compensation payments associated with deviations from ideal service times or pickup at collection points. The hierarchical nature of the interaction between the company and its customers makes a bilevel optimization framework a natural modeling approach for the \ABC. In summary, the \ABC\ integrates routing and facility location decisions while accounting for customer preferences and compensation mechanisms.\\

The main contributions of this paper are summarized as follows:
\begin{itemize}
    \item We introduce the \ABC, a bilevel location-routing problem that integrates collection-point activation, routing, customer choice for service mode, and compensation mechanisms for service inconvenience.

    \item We exploit the structure of the customers' problems to derive two equivalent single-level MILP reformulations. The first one is obtained from primal/dual optimality conditions, whereas the second one uses primal variables only and incorporates a system of inequalities that captures the customers preferences.

    \item We develop a two-phase MIP heuristic and a warm-start strategy to improve the solution of larger instances.

    \item We conduct extensive computational experiments and quantify the impact of explicitly incorporating customer preferences. 
    
    \item A comparison of the \ABC\ with a benchmark model that ignores customers' preferences indicates that explicitly accounting for customer choice alters the resulting location-routing configuration while improving customer preference satisfaction, with only a small increase in the company's overall cost.
\end{itemize}

The remainder of this paper is structured as follows. Section \ref{sec:lit_rev} reviews the literature on contributions related to the \ABC. Section 3 defines the problem and discusses its main features through an illustrative example, while Section \ref{sec:Bilevel} formulates the problems at both decision levels and presents the proposed bilevel mathematical model. Section \ref{sect:single-level} introduces an equivalent non-linear single-level reformulation as well as the two alternative linearizations. The two-phase MIP heuristic is described in Section  \ref{sec:heur}. Section \ref{sect:results} summarizes the computational experiments that have been conducted and analyzes the results obtained. Managerial insights are discussed in Section \ref{sect:managerial}: The comparison of \ABC\, against a benchmark model that ignores customers preferences is presented in Section \ref{sect:ValoCust}, whereas in Section \ref{sect:sol-analysis} the optimal solution to a small case study with 40 customers and 8 collection points is used to illustrate the trade-off between routing efficiency and customers' satisfaction. The manuscript ends with Section \ref{sect:conclu}, where the conclusions and future research directions are outlined.

\section{Literature Review} \label{sec:lit_rev}

Last-mile delivery systems increasingly combine home delivery with out-of-home alternatives, such as parcel lockers and collection-and-delivery points (CDPs). These alternatives enable the consolidation of multiple deliveries at common locations, potentially reducing the routing effort associated with traditional home delivery. Consequently, a growing body of research has investigated how home delivery and pickup options can be jointly incorporated into last-mile distribution planning.\\

Several routing problems have incorporated this combination of home delivery and alternative pickup locations. For example, \citet{grabenschweiger2021vehicle} study a vehicle routing problem with heterogeneous locker boxes in which customers can be served either at home within a time window or through a locker station. Similarly, \citet{tilk2021last} introduce the vehicle routing problem with delivery options, where customers may be associated with alternative delivery locations and time windows and the carrier selects the final delivery option while satisfying a prescribed service level regarding customer preferences. \citet{dumez2021large} also consider multiple delivery options characterized by preference levels and time windows, while accounting for capacities at shared delivery locations. \citet{janinhoff2024multitrip} consider a multitrip setting combining home delivery, parcel lockers, and pickup stations, while \citet{dell2023pickup} integrate parcel lockers into a pickup-and-delivery problem with time windows. These studies show how alternative delivery locations can provide additional flexibility to last-mile routing. However, the locations of the pickup facilities are generally treated as given rather than being jointly determined with the delivery routes.\\

While the previous studies assume that the locations of the alternative pickup facilities are given, other studies have integrated location into last-mile distribution planning. In this direction, \citet{janjevic2019integrating} incorporate the location of collection-and-delivery points within the strategic design of an e-commerce distribution network, considering the interaction between facility placement and the resulting distribution operations. More closely related to a location-routing setting, \citet{grabenschweiger2022multi} consider a multi-period problem in which locker locations are selected while home-delivery, locker-assignment, and routing decisions are jointly determined. More recently, \citet{ozyavas2025designing} study a delivery network in which parcel lockers can serve both as customer collection points and as transfer points for the delivery company, jointly considering their location and the associated routing decisions. These studies highlight the importance of integrating the design of the pickup network with last-mile distribution decisions rather than considering facility location and vehicle routing separately.\\

Related problems have also been studied within two-echelon distribution systems. For instance, \citet{zhou2018multi} consider a multi-depot two-echelon vehicle routing problem with delivery options, in which goods are transported from depots to satellites and subsequently delivered to customers, while intermediate pickup facilities provide an alternative for parcel collection. Similarly, \citet{enthoven2020two} integrate parcel lockers and cargo bikes into a two-echelon distribution system, where customers can either collect their parcels from covering locations or receive home delivery through vehicles operating from satellite facilities. Although these problems share some characteristics with the \ABC, their distribution structure is fundamentally different. In two-echelon systems, intermediate facilities primarily support the consolidation or transfer of goods between distribution stages. In the \ABC, collection points instead constitute alternative final service locations, and the same route can visit both activated collection points and customers offered home delivery.\\

Beyond the configuration of the delivery network, customer preferences constitute another important dimension of last-mile delivery. This is particularly relevant for attended home delivery, where customers may value delivery times differently and the attractiveness of home delivery therefore depends on when the service is provided. \citet{agatz2011time} highlight the trade-off between operational efficiency and customer service when designing delivery time slots. More recently, \citet{amorim2024customer}, using empirical data from an online grocery retailer, show substantial heterogeneity in customer preferences for delivery-service attributes. In particular, they find that delivery speed, time-slot precision, and timing are important drivers of customers' time-slot selections. These findings support the explicit consideration of time-related customer preferences in last-mile delivery decisions rather than treating feasible delivery times as equally attractive to customers.\\

From an optimization perspective, customer preferences have been incorporated directly into last-mile delivery models. \citet{tilk2021last} allow customers to express preferences among alternative delivery locations and time windows; however, the carrier ultimately selects the delivery option while satisfying a prescribed service level regarding customer preferences. In attended home delivery, \citet{mackert2019integrating} explicitly incorporate customer choice among delivery time slots, recognizing that the alternatives offered by the retailer influence the subsequent customer response. More recently, \citet{zhou2026last} consider both home-delivery time slots and out-of-home locations under heterogeneous customer preferences. In their setting, the retailer determines which delivery alternatives and prices to offer, after which customers select among the available options according to a stochastic choice model. These studies represent an important shift from treating customer preferences only as service attributes toward explicitly accounting for how customers respond to the delivery alternatives offered by the logistics provider. \\

Customer preferences are often associated with compensation mechanisms offered to customers. For instance, \citet{grabenschweiger2021vehicle} and \citet{grabenschweiger2022multi} include compensation costs associated with serving customers through locker boxes, while \citet{dell2023pickup} consider discounts granted to customers choosing self-service at lockers. Similarly, \citet{enthoven2020two} incorporate customer-specific connection costs reflecting the attractiveness of alternative delivery modes. More recently, \citet{peng2026distributionally} explicitly consider financial incentives offered to customers to encourage delivery postponement or switching from home delivery to pickup alternatives. These approaches internalize customer-related service consequences in the objective function, although the associated costs are typically linked to the selected delivery mode, service adjustment, or incentive offered rather than to explicit ranked preferences over both delivery times and collection points.\\

When customers retain autonomy over their final service choice, the interaction between the service provider and the customers can naturally be represented as a hierarchical decision process. Bilevel optimization provides a suitable framework for such settings, as the provider first determines the alternatives made available and customers subsequently respond according to their own preferences. This perspective has been explored in the facility-location literature. \citet{camacho2014solving} consider an uncapacitated facility location problem in which the company determines which facilities to open, while customers are assigned according to their ordered preferences over the available facilities. \citet{casas2018approximating} extend this setting to capacitated facilities, with customers seeking service from their most preferred facilities according to predetermined preference lists. Similarly, \citet{calvete2020matheuristic} incorporate restrictions on the number of customers that can be allocated to each facility while preserving customer autonomy in selecting their most preferred open facility. These studies demonstrate how bilevel optimization can explicitly capture hierarchical interactions between a service provider and customers whose decisions are driven by ranked preferences.\\

Bilevel optimization has also been used to model hierarchical interactions in last-mile delivery. \citet{cerulli2024bilevel} consider a peer-to-peer logistics platform in which the platform assigns delivery requests to independent carriers, who subsequently decide which requests to accept according to their own profitability. More closely related to customer choice, \citet{peng2026distributionally} study a last-mile delivery system in which the operator offers financial incentives to encourage customers to accept delayed delivery or switch from home delivery to alternative pickup modes. Customers then decide whether to accept these offers, while the operator anticipates their responses when planning the delivery operations. Although the followers and their decision mechanisms differ across these applications, these studies illustrate the relevance of bilevel optimization when last-mile decisions depend on the autonomous responses of other actors.\\

To synthesize the literature discussed above, Table~\ref{tab:lit_comp} compares the studies most closely related to this work in terms of their main characteristics, treatment of customer preferences, and representation of customer response.\\

\begin{table}[htbp]
\centering
\caption{Comparison of related studies in last-mile delivery.}
\label{tab:lit_comp}
\resizebox{\textwidth}{!}{
\begin{tabular}{lccccccccl}
\hline
\textbf{Reference} &
\textbf{HD} &
\textbf{Pickup} &
\textbf{CP } &
\textbf{Routing} &
\textbf{Time } &
\textbf{CP } &
\textbf{Preference-} &
\textbf{Customer} &
\textbf{Choice } \\

& & & \textbf{location} & & \textbf{ranking} &
\textbf{ranking} &
\textbf{related cost} &
\textbf{response} &
\textbf{mechanism} \\

\hline

\citet{zhou2018multi}
& $\checkmark$ & $\checkmark$ & -- & $\checkmark$
& -- & -- & -- & --
& Centralized selection \\

\citet{enthoven2020two}
& $\checkmark$ & $\checkmark$ & -- & $\checkmark$
& -- & -- & $\checkmark$ & --
& Centralized selection \\

\citet{grabenschweiger2021vehicle}
& $\checkmark$ & $\checkmark$ & -- & $\checkmark$
& -- & -- & $\checkmark$ & --
& Centralized selection \\

\citet{tilk2021last}
& $\checkmark$ & $\checkmark$ & -- & $\checkmark$
& -- & -- & -- & --
& Priority-based service level \\

\citet{grabenschweiger2022multi}
& $\checkmark$ & $\checkmark$ & $\checkmark$ & $\checkmark$
& -- & -- & $\checkmark$ & --
& Centralized selection \\

\citet{dell2023pickup}
& $\checkmark$ & $\checkmark$ & -- & $\checkmark$
& -- & -- & $\checkmark$ & $\checkmark^{a}$
& Customer/company selection \\

\citet{janinhoff2024multitrip}
& $\checkmark$ & $\checkmark$ & -- & $\checkmark$
& -- & -- & -- & $\checkmark^{b}$
& Product-based eligibility \\

\citet{ozyavas2025designing}
& -- & $\checkmark$ & $\checkmark$ & $\checkmark$
& -- & -- & -- & --
& Centralized selection \\

\citet{zhou2026last}
& $\checkmark$ & $\checkmark$ & -- & $\checkmark$
& -- & -- & -- & $\checkmark$
& Stochastic utility maximization \\

\citet{peng2026distributionally}
& $\checkmark$ & $\checkmark$ & -- & $\checkmark$
& -- & -- & $\checkmark$ & $\checkmark$
& Incentive acceptance \\

\hline
\textbf{This study}
& $\checkmark$ & $\checkmark$ & $\checkmark$ & $\checkmark$
& $\checkmark$ & $\checkmark$ & $\checkmark$ & $\checkmark$
& Ranked preferences \\
\hline
\end{tabular}
}

\begin{flushleft}
\footnotesize
HD: home delivery; CP: collection point.
Time ranking and CP ranking indicate explicit ordered preference lists over
home-delivery time intervals and collection points, respectively.
Preference-related cost indicates whether customer preferences or service
inconvenience are explicitly represented in the objective function through
compensation, discounts, incentives, or other customer-specific cost terms.
Customer response indicates that at least part of the final service decision
is made by the customer rather than entirely by the logistics provider.
$^{a}$ Customers may select home service or self-service at a nearby locker,
or delegate the service decision to the logistics company.
$^{b}$ Customers select a delivery product that determines the set of eligible
delivery locations; when multiple locations are eligible, the final location
may be selected by the service provider.
\end{flushleft}
\end{table}

The comparison in Table~\ref{tab:lit_comp} shows that, although previous studies have addressed several components of the \ABC, important differences remain in how customer preferences, preference-related costs, and customer responses are integrated with operational decisions. Location-routing models typically retain control over the final delivery alternative, whereas models with endogenous customer choice generally assume a predetermined pickup infrastructure. Moreover, customer preferences are commonly represented through service levels, eligibility rules, stochastic utility, or incentive mechanisms, and preference-related costs are typically associated with compensation, discounts, or service adjustments rather than with explicit rankings over both home-delivery time intervals and collection points.\\

The \ABC\ bridges these perspectives by jointly considering collection-point location, routing and scheduling, and endogenous customer responses. Customers rank jointly home-delivery time intervals and candidate collection points, and select among the alternatives made available by the company accordingly. This interaction naturally leads to a hierarchical decision process in which the company anticipates customer choices when planning its operations.

\section{Problem Definition}\label{sec:Prob-def}

Consider a graph $G=(V, A)$ with $V= I \cup K \cup \{0 \}$, where $I$ denotes the set of customers with service demands (parcels), $K$ the set of potential collection points, and $0$ a central depot. The arc set is defined as $A = \{ (i,j) : i, j \in V, i \ne j \}$.\\

Each arc $(i,j) \in A$ has a routing cost $c_{ij}$ and a travel time $t_{ij}$. Furthermore, each potential collection point $k \in K$ has an activation (setup) cost $f_k$.\\

Customers can be served either at their respective (home) locations or at an activated collection point. Each customer $i \in I$ has an ideal service time, $T^*_i$, for receiving their parcel at home, but may accept home delivery at a different time according to their preferences. 
Specifically, the time horizon is partitioned into intervals identified by customer $i \in I$ as $\Delta_{ir} = (\delta_{i,\,r-1}, \delta_{ir}]$, indexed by the set $R^i = \{0, 1, \dots, \bar{r}^i\}$. 
Note that each customer $i \in I$ may have a different number of time intervals during which they are willing to receive their parcel ($\bar{r}^i$). A service time falling within interval $\Delta_{ir}$ is associated with a preference value $p^T_{ir}$ for being served during interval $r \in R^i$. For customer $i \in I$, if $p^T_{ir} < p^T_{ir'}$, this indicates that interval $r$ is preferred over interval $r'$. {In general, these preferences are represented by integer values, where 0 corresponds to the interval containing the ideal delivery time and therefore represents the most preferred interval. Similarly, $\bar{r}^i$ corresponds to the least preferred interval for home delivery.}   \\

{Additionally, each customer $i \in I$ may patronize collection points according to their own preferences. Consequently, customers may prefer to pick up their parcels at an activated collection point rather than receive them at home.} To capture this flexibility in the delivery system, we define $p^F_{ik}$ as the preference of customer $i \in I$ for collection point $k \in K$. Similarly to time-related preferences, each customer $i \in I$ ranks all collection points according to their preferences, where 1 denotes the most preferred collection point and $|K|$ the least preferred. Hence, $p^F_{ik} < p^F_{ik'}$ if and only if customer $i \in I$ prefers collection point $k \in K$ over $k' \in K$.\\ 


As mentioned above, when a parcel is not delivered to the customer's location at their ideal service time, the company pays compensation to the customer. Compensations are determined based on the company's estimate of the inconvenience caused to customers by the service provided, using the customer information available to the company. We assume that the company knows only the customers' locations and ideal service times and does not have access to their preference lists over alternative time intervals or collection points. Thus, for customers ultimately served at home, compensation depends on the deviation from their ideal service time, whereas for customers served through a collection point, compensation depends on the distance to the closest activated collection point. \\

Specifically, we denote by $\theta^T_i$ the compensation paid to customer $i\in I$ per unit of deviation from the ideal service time when the parcel is delivered to the customer's location. For customers served through a collection point, the compensation is defined as
$\min\{\theta_{ik}^F: k\in K \textit{ is an activated collection point}\}$, where $\theta^F_{ik}$ is a parameter directly proportional to the distance between the customer's location and collection point $k$. This rule applies both to customers who are not offered home delivery and to those who reject the proposed home-delivery service and choose pickup instead. In either case, the customer selects an activated collection point according to their preferences, while the compensation is based on the distance to the closest activated collection point, using only the information available to the company.
\\

Given that compensation is based solely on the information available to the company, namely the customer's location and ideal service time, it is independent of other customers' preferences.\\



Another crucial aspect is that a customer $i \in I$ may prefer certain collection points to certain service-time intervals. Therefore, for each customer $i \in I$, we consider a merged preference list in which time intervals $r \in R^i$ and collection points $k \in K$ are jointly ranked according to a relative preference order over potential \textit{service modes} (i.e., alternative delivery time intervals or collection points). In this merged list, defined as $\gamma_{il}$ with $\ell \in L^i=\{0, 1,2,\dots,|R^i| + |K|-1\}$, the components associated with the first $|R^i|$ indices indicate preferences for service time deviations in the $l$-th interval, whereas the components associated with indices  $\ell \in\{|R^i|, \dots, |R^i|+|K|-1\}$ indicate preferences for collection point  $k=\ell-|R^i|+1$. Relative preferences take (integer) values ranging from 0 (the most preferred service mode) to $|R^i|+|K|-1$ (the least preferred service mode). Indeed, the best-ranked service mode for each customer is to receive the parcel at home in their ideal time, so $\gamma_{i0} = 0$ and their least-ranked service mode is to receive the parcel after their time window availability $\gamma_{i,|R^i|-1} = |R^i|+|K|-1$. 
 The remaining service modes are ordered based on the customer’s relative preferences, so they may indistinctly represent deviated time intervals or  collection points. In the end, we obtain a merged list of preferences such that $\gamma_{i \ell} <\gamma_{i \ell'}$, if and only if customer $i \in I$ prefers service mode $\ell$ to service mode $\ell'$, where $\ell, \ell'\in L^i$ and each of them may represent either  a service time deviation in the interval $\ell$ (when $0\leq \ell \leq |R^i|-1$) or collection point $\ell-|R^i|+1$  (when $\ell> |R^i|$); and similarly for $\ell'$.\\

In the \ABC\ the company must determine the set of collection points to activate, as well as the set of customers who may potentially receive home delivery. In addition, the company must design a service route that starts and ends at the depot and visits all activated collection points as well as all customers selected for home delivery.\\

Customers who are not selected for home delivery must collect their parcels from an activated collection point according to their preferences. Moreover, a customer may reject home delivery if the offered delivery time is less preferred than another available option. That is, if there exists an activated collection point at which the customer prefers to pick up the parcel rather than receive it at home at the proposed delivery time, the customer will choose that collection point instead. In such a case, the customer will select their most preferred activated collection point.\\

Customers who ultimately receive home delivery receive compensation based on the deviation of their delivery time from their ideal service time, whereas  customers who receive service at collection points, receive compensation based on the closest activated collection point.\\

\noindent The decisions of the company are the following:
\begin{itemize}
    \item[-] A set of collection points to be activated.
    \item[-] A set of customers selected to potentially receive home delivery. 
   \item[-] A service route, starting and ending at the depot, that visits all selected customers as well as all activated collection points.
    \item[-]  A home-delivery time for each selected customer, as determined by the service route.
\end{itemize}

\noindent Given a company’s solution, each customer must decide:
\begin{itemize}
    \item[-] Whether to accept the proposed home-delivery time (if home delivery is offered).
    \item[-] Which activated collection point to use for parcel pickup when home delivery is rejected or not offered.
\end{itemize}

\noindent The objective function minimizes the following cost-related terms:
\begin{itemize}
\item[-] The setup cost of the activated collection points.
\item[-] The total routing cost for the selected customers and the activated collection points.
\item[-] The economic compensations paid to customers for the service they receive.
\end{itemize}

The \ABC\ consists of determining a company's solution that minimizes the cost objective function described above while accounting for the customers' optimal reactions based on their preferences.

\subsection{Illustrative Example}\label{Ex1-bis}
\begin{figure}[H]
    \centering
    \captionsetup{font=small, labelfont=bf, skip=6pt}
    \captionsetup[subfigure]{
        font=small,
        labelfont=bf,
        justification=centering,
        singlelinecheck=false,
        skip=4pt
    }

    \begin{subfigure}[t]{0.5\textwidth}
        \centering
        \includegraphics[width=0.86\linewidth]{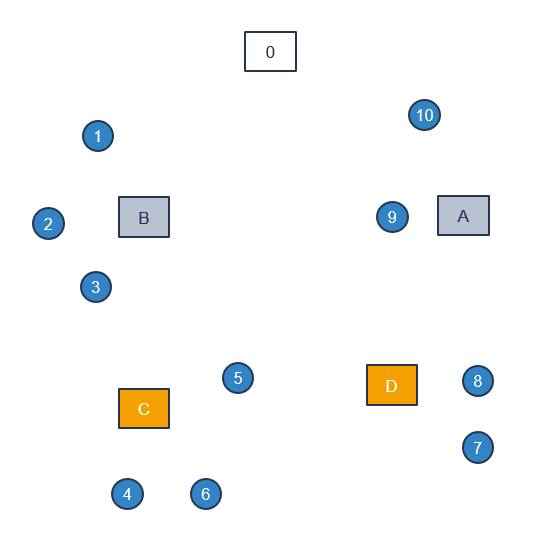}
        \caption{Input graph}
        \label{InputGraph-bis}
    \end{subfigure}%
    \begin{subfigure}[t]{0.5\textwidth}
        \centering
        \includegraphics[width=0.93\linewidth]{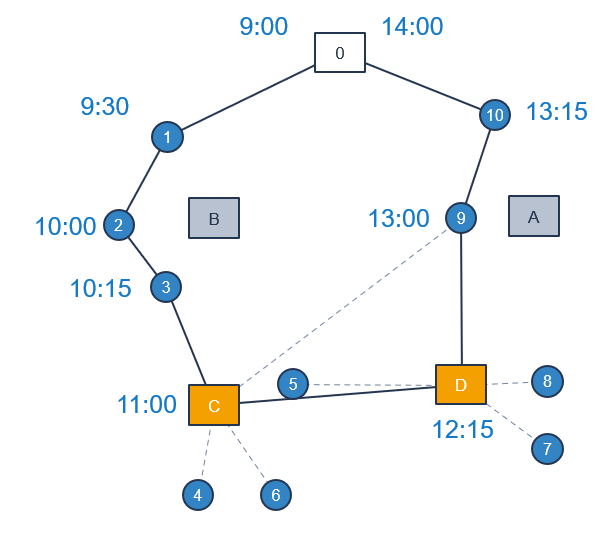}
        \caption{Leader's route through activated facilities and customers accepting the scheduled service time}
        \label{Solution-bis}
    \end{subfigure}

    \caption{Figure for Example~\ref{Ex1-bis}}
    \label{Fig1}
\end{figure}

Consider the example illustrated in Figure \ref{InputGraph-bis}, with ten customers indexed in $I = \{1, 2, 3, 4, 5, 6, 7, 8, 9, 10\}$, four potential collection points in $K = \{A, B, C, D\}$, and one depot located at node $0$.
We assume that the company offers service in the time interval 9h-13h, and $f_k=0$ for all $k\in K$.  
Time-related customer preferences ($p^T_{ir}$) are given in Table \ref{time-compensation-tris}, whereas {Table \ref{facility-compensation-bis}} provides information on preferences related to collection points ($p^F_{ik}$). \\


In Table \ref{time-compensation-tris}, each row corresponds to a customer and contains at least four entries: the first indicates its ideal service time ($T^*_i$), and the second indicates the number of admissible time intervals ($\bar r^i-1$). The remaining entries report the $\bar r^i-1$ admissible time intervals, ordered chronologically, along with their respective preferences. For example, customer 3 has an ideal service time at 10:00h and two admissible time intervals. The first interval, [8:00, 9:00], has a preference rank of 2, whereas the second interval, [10:00, 14:00], has a preference rank of 1. This indicates that the customer prefers delivery within the latter interval, which includes their ideal service time, but is also willing to accept an earlier delivery within their second-most-preferred interval. \\

\begin{table}[H]
	\centering
	\begin{tabular}{c|c|c|ll|cl}
	Customer $(i)$ & $T_i^*$ (h) & $\bar r^i-1$ & \multicolumn{1}{c}{$\Delta_{i1}$} & {$p^T_{i1}$} & {$\Delta_{i2}$} & {$p^T_{i2}$} \\
	\hline
	1  & 10:00 & 2 & [8:00, 10:00] & 1   & [10:00, 11:00] & 2\\
	2  & 12:00 & 1 & [9:00, 12:00] & 1   &  --\\
	3  & 10:00 & 2 & {[8:00, 9:00]}& {2} & {[10:00, 14:00]}  & {1}\\
	4  & 10:00 & 1 & [8:00, 10:00]  & 1 & -- \\
	5  & 9:00  & 2 & {[8:00, 9:00]} & {2}& {[15:00, 16:00]} & {1}\\
	6  & 12:00 & 2 & [8:00, 12:00]   & 1 & [15:00, 16:00] & 2\\
	7  & 10:00 & 1 & [10:00, 11:00]  & 1 & -- &  \\
	8  & 9:00  & 1 & [8:00, 9:00] & 1 &  -- &  \\
	9  & 12:00 & 2 & [11:00, 12:00]  & 1& [12:00, 16:00] & 2\\
	10 & 14:00 & 1 & [8:00, 14:00]  & 1&  -- &  \\
	\end{tabular}
	\caption{Service time windows for customers. The ideal service time $T_i^*$ (in hours) and the number of admissible time intervals $\bar r^i-1$ are reported for each customer. The remaining entries list the admissible time windows in order of preference.}
	\label{time-compensation-tris}
    
\end{table}

Analogously, Table \ref{facility-compensation-bis} shows the assignment preferences $p^F_{ik}$ of each customer over the available collection points. Each row represents a customer, and the entries indicate the corresponding preference ranking, where 1 denotes the most preferred collection point and higher values indicate lower preference levels. For example, customer 9 prefers collection point A, followed by collection points C, D, and B. \\

Table \ref{Merged2-bis-bis} shows the merged lists of relative preferences. For each customer $i \in I$, the column labeled $|L^i|$ gives the number of elements in the corresponding merged list. The following block of columns reports the joint preference rankings assigned to the different service modes (home delivery or collection-point pickup). For each customer, the entries in the first $|R^i|$ columns of this block (before the double vertical line), corresponding to indices $l \in R^i$, represent preferences for service time deviations within interval $\Delta_{il}$. The last $|K|$ entries (after the double vertical line), corresponding to indices $|R^i| < l \leq |L^i|$, represent preferences for collection points, where $k = l - |R^i| \in K$. As can be seen, for all customers, the most preferred option is to be served exactly at the ideal service time ($l=0$ with $\gamma_{i0}=0$), whereas the least preferred option is to be served outside the set of acceptable time windows ($l=|R^i|$ with $\gamma_{i|R^i|}=|L^i|$). The ranking of the remaining options varies across time-related service modes and collection-point service modes. Note that this ranking does not necessarily alternate between the two types of service modes and depends on each customer's preferences.\\

\begin{table}
\begin{subtable}{.5\linewidth}
	\centering
	\begin{tabular}{c|cccc}
	Customer $(i)$ & A & B & C & D \\
	\hline
	1  & 1 & 3 & 2 & 4 \\
	2  & 1 & 2 & 3 & 4 \\
	3  & 4 & 2 & 1 & 3 \\
	4  & 2 & 4 & 1 & 3 \\
	5  & 3 & 4 & 2 & 1 \\
	6  & 2 & 1 & 3 & 4 \\
	7  & 3 & 4 & 2 & 1 \\
	8  & 1 & 2 & 4 & 3 \\
	9  & 1 & 4 & 2 & 3 \\
	10 & 3 & 4 & 1 & 2 \\
	\end{tabular}
	\caption{Collection-point preferences $p^F_{ik}${, $k\in K$}  for customers.\\ For each customer, the available collection points are ranked according to preference, where 1 denotes the most preferred option.}
	\label{facility-compensation-bis}
\end{subtable}
\hspace{1.cm}\begin{subtable}{.5\linewidth}
		\centering
		\begin{tabular}{c|c|ccccccc}
        &&\multicolumn{7}{c}{Service mode ($l$)}\\
$i$& $|L^i|$ & 0 & 1 & 2 & 3 & 4 & 5 & 6\\
\hline
1 & 7 & $0$ & $1$ & $6$ & \multicolumn{1}{||c}{$2$} & $4$ & $3$ &  {$5$}\\
2 & 6 & $0$ & $2$ & \multicolumn{1}{||c}{$1$} & $3$ & $4$ &  {$5$}\\
3 & 7 & $0$ & {$6$} &{$1$}&  \multicolumn{1}{||c}{$5$} & $3$ & $2$ &  $4$\\
4 & 6 & $0$ & $1$ &\multicolumn{1}{||c}{$3$} & $5$ & $2$&  {$4$}\\
5 & 7 & $0$ & {$6$} & {$3$} & \multicolumn{1}{||c}{$4$} & $5$ & $2$ &  {$1$} \\
6 & 7 & $0$ & $1$ &  $6$ & \multicolumn{1}{||c}{$3$} & $2$ & $4$ &  {$5$}\\
7 & 6 & $0$ & $5$ & \multicolumn{1}{||c}{$3$} & $4$ & $2$ & {$1$}\\
8 & 6 & $0$ & $1$ & $6$ & \multicolumn{1}{||c}{$2$} & $3$ & $5$ &  {$4$}\\
9 & 7 & $0$ & $3$ & $6$ & \multicolumn{1}{||c}{$1$} & $5$ & $2$ &  {$4$}\\
10& 6 & $0$ & $1$ & \multicolumn{1}{||c}{$4$} & $5$ & $2$ &  {$3$}\\
\end{tabular}
\caption{Merged preference lists {$\gamma_{i}^l$, $l\in L^i$} for customers.\\$\phantom{asdaafs}$\\{\color{white}asasad}}
\label{Merged2-bis-bis}
\end{subtable}
\end{table}

Figure \ref{Solution-bis} shows a possible solution for the company. In this solution, collection points $C$ and $D$ are activated, and home delivery is offered to customers ${1, 2, 3, 9, 10}$. The numbers next to the nodes indicate the arrival times of the route at the corresponding locations. The route departs from the depot at $9\text{:}00$h. The arrival times at the selected customers are as follows. For customer 1, the arrival time is $9\text{:}30$ h, resulting in a deviation of $h_1=|9.5-10|=0.5$ hours from the ideal delivery time. Similarly, customers 2, 3, 9, and 10 have arrival times of $10\text{:}00$, $10\text{:}15$, $13\text{:}00$, and $13\text{:}15$ h, respectively. Their corresponding deviations from the ideal delivery times are $h_2=|10-12|=2$ hours, $h_3=|10.25-10|=0.25$ hours, $h_9=|13-12|=1$ hour, and $h_{10}=|13.25-14|=0.75$ hours. Note that arrival times are reported in clock time and denoted by $g_i$ when expressed in decimal units. Deviations from the ideal delivery time $T^*_i$ are calculated as $h_i=|g_i-T^*_i|$ and are also expressed in decimal units.\\  
  
The customers respond according to their preferences. Customers 1, 2, 3, and 10 accept the company's offer because $g_1=9.5 \in {\Delta_{11}}$, $g_2=10 \in {\Delta_{21}}$, $g_3=10.25 \in {\Delta_{32}}$, and $g_{10}=13.25 \in {\Delta_{{10},\,1}}$. In all four cases, the proposed service-time interval is the most preferred one, $p^T_{11}=p^T_{21}=p^T_{32}=p^T_{10,1}=1$, and it is preferred to picking up the parcel at either of the activated collection points (C or D).
For instance, note that customer 2 prefers picking up the parcel at collection point A to receiving it at home at 10:00h $\left(2=\gamma_{21}>\gamma_{22}=1\right)$, but A is not activated. Therefore, customer 2 selects the second-most preferred option. In contrast, customer 9 rejects the offered service time because $g_9=13 \in {\Delta_{92}}$, which is less preferred than picking up the parcel at the activated collection point C $\left(6=\gamma_{92}>\gamma_{25}=2\right)$. However, the customer receives the compensation associated with collection point D, i.e., $\theta^F_{9D}$, since D is closer to customer 9 than C. Recall that customers not served through home delivery are compensated based on the closest activated collection point, even though they pick up their parcels at their most preferred activated collection point. For example, customer 5 prefers collection point D over C but is closer to C and therefore receives $\theta^F_{5C}$. For the remaining customers, the closest activated collection point coincides with the most preferred one. Table~\ref{tab:summary} summarizes the company's offers and the customers' responses. As can be seen, the overall compensation paid by the company is $\theta^T_1h_1+\theta^T_2h_2+\theta^T_3h_3 +  \theta^F_{4C} + \theta^F_{5C} + \theta^F_{6C} + \theta^F_{7D}+\theta^F_{8D} + \theta^F_{9D} +\theta^T_{10}h_{10}$.

\begin{table}[H]
	\centering
	\begin{tabular}{c|cc|ccc|cc}
	Customer & & & \multicolumn{3}{c|}{Preferences with respect to:} & Customers' & Compensation\\
$i$& Arrival time & $h_i$ & \multicolumn{1}{c}{Arrival time} & \multicolumn{1}{c}{C} & \multicolumn{1}{c|}{D} & \multicolumn{1}{c}{Response} & \multicolumn{1}{c}{Paid}\\
\hline
1 & 9:30  & $0.5$    & $1$ & $2$ & 5 &  Accept & $0.5\theta^T_1$\\
2 & 10:00 & $2$      & $2$ & $4$ & 5 & Accept & $2\theta^T_2$\\
3 & 10:15 & $0.75$   & $1$ & $2$ & 4 & Accept & $0.75\theta^T_3$\\
4 & --    & --       & -- & 2 & 4 & C & $\theta^F_{4C}$\\
5 & --    & --       & -- & 2 & 1 & D & $\theta^F_{5C}$\\ 
6 & --    & --       & -- & 4 & 5 & C & $\theta^F_{6C}$\\
7 & --    & --       & -- & 2 & 1 & D & $\theta^F_{7D}$\\
8 & --    & --       & -- & 5 & 3 & D & $\theta^F_{8D}$\\
9 & 13:00 & $1$      & $6$ & $2$ & 4 & C & $\theta^F_{9D}$\\
10& 13:15 & $0.75$   & $1$ & $2$ & 3 & Accept & $0.75\theta^T_{10}$\\
\end{tabular}
\caption{Customers' responses to the company's proposed service options.}
\label{tab:summary}
\end{table}

\section{A Bilevel Programming Formulation for the LRP-DP}\label{sec:Bilevel}

As mentioned above, the \ABC clearly exhibits a hierarchical relationship among decision-makers. Therefore, it naturally fits within a bilevel optimization framework, where the company (hereafter referred to as the leader) is responsible for the upper-level decisions, while the customers (hereafter referred to as the followers) make the lower-level decisions. Since customers are independent and do not communicate with each other, they can be aggregated into a single problem encompassing all of them \citep{calvete2007linear}.\\

In addition to the notation introduced in Section \ref{sec:Prob-def}, we define the following: for a given node set $S \subseteq V\setminus\{0\}$, $\delta^+(S) = \{(i, j) \in A : i \in S, j \notin S\}$ and $\delta^-(S) = \{(i, j) \in A : i \notin S, j \in S\}$ denote the sets of arcs \emph{leaving} and \emph{entering} the node set $S$, respectively. In particular, when $S$ is a singleton, i.e., $S = \{i\}$ with $i \in V$, we write $\delta^+(i) = \{(i,j) \in A\}$ and $\delta^-(i) = \{(j,i) \in A\}$ to denote the sets of arcs \emph{leaving} and \emph{entering} node $i \in V$, respectively.\\

The decision variables related to the leader are the following:
\begin{itemize}[label={}]
\item $y_{k}=\begin{cases}1& \text{if collection point $k\in K$ is activated (so it must be visited by the route)} \\ 
0 & \text{otherwise.}
\end{cases}$\\

\item $s_{i}=\begin{cases}1& \text{if customer $i\in I$ is visited by the service route} \\ 
0 & \text{otherwise.}
\end{cases}$\\

\item $x_{ij}=\begin{cases}1& \text{if the service route  uses arc } (i,j)\in A \\ 
0 & \text{otherwise.}
\end{cases}$\\

\item[] $g_i$: time when the service route arrives to node $i\in V$. For non-visited nodes, this variable may take any value as it is irrelevant.
For the central depot, we set $g_0=0$, as we assume the route departs from the depot at time 0.
\item[] $h_i$: deviation (in absolute value) of the arrival time to customer $i \in I$ from its ideal service time $T^*_i$; that is, $h_i = |g_i - T^*_i|$.\\ 
\end{itemize}

Denote a leader's decision as {$\bm\ell=\left(\bm y,\bm s, \bm x, \bm g, \bm h\right)=\left(y_k,s_i,x_{ij},g_i,h_i\right)_{i\in I, (i,j)\in A, k\in K}$}. Given the parameters introduced previously and the recently defined decision variables, the feasible solutions for the leader must satisfy the following constraints:
\begin{subequations}
\begin{align}
 & x(\delta^+(0))=  x(\delta^-(0))= 1 && \label{eq:L_depot-0}\\
&x(\delta^+(i))= s_i\qquad && i\in I \label{eq:L_llegar_cust-0}\\
 & x(\delta^-(i))= s_i\qquad && i\in I \label{eq:L_salir_cust-0}\\
& x(\delta^+(i))= y_k\qquad && k\in K \label{eq:L_llegar_collec-0}\\
 & x(\delta^-(i))= y_k\qquad && k\in K \label{eq:L_salir_collec-0}\\
 & g_j\geq g_i + t_{ij}- \text{M}(1-x_{ij})\qquad && (i, j)\in A, j\ne 0 \label{eq:L_arrival_time-0}\\
  &\delta_{i0}- g_i \leq \text{M} (1- s_i) \qquad && i\in I    \label{eq:F_antesTmin-0}\\
 & g_i-\delta_{i\bar r^i}\leq \text{M} (1- s_i)\qquad && i\in I    \label{eq:F_e_despuesTmax-0}\\ 
  & h_i \geq  g_i - T^*_i \qquad && i\in I     \label{eq:F_h_izq-0}\\ 
 & h_i \geq T^*_i -  g_i. \qquad && i\in I    \label{eq:F_h_der-0} 
 \end{align}
\end{subequations} 

 Constraints \eqref{eq:L_depot-0} ensure that the route starts and ends at the depot. Constraints \eqref{eq:L_llegar_cust-0} and \eqref{eq:L_salir_cust-0} ensure that, for each selected customer, exactly one arc enters the corresponding node and exactly one arc leaves it. Constraints \eqref{eq:L_llegar_collec-0} and \eqref{eq:L_salir_collec-0} operate analogously for the activated collection points. The Miller--Tucker--Zemlin (MTZ) constraints \eqref{eq:L_arrival_time-0}, where M is a sufficiently large positive constant, regulate the arrival times at the visited customers and prevent subtours. If customers have earliest and latest allowable delivery times, these requirements are enforced by constraints \eqref{eq:F_antesTmin-0} and \eqref{eq:F_e_despuesTmax-0}, respectively. In our case, $T^{\min}_i=\delta_{i0}$ and $T^{\max}_i=\delta_{i\bar r^i}$. Finally, constraints \eqref{eq:F_h_izq-0} and \eqref{eq:F_h_der-0} define the deviations with respect to the ideal service times. In the following, we denote by $\mathcal{X}^{\bm\ell}$ the feasible domain determined by  solutions $\bm\ell=\left(\mathbf y,\mathbf s, \mathbf x, \mathbf g, \mathbf h\right)$ with  $\mathbf y,\, \mathbf s,\, \mathbf x$, binary and $\mathbf g\geq \mathbf 0, \mathbf h\geq \mathbf 0$ that satisfy constraints \eqref{eq:L_depot-0}- \eqref{eq:F_h_der-0}.\\

The above variables allow modeling the portion of the objective function that depends on the leader's decisions, namely facility activation costs and routing costs. However, the compensation costs cannot be modeled unless the customers' responses are known. 
In particular, determining whether customer $i \in I$ receives home delivery or selects a collection point is crucial, as this determines the compensation paid by the leader to customer $i$.
Note that this information depends not only on the offered service time but also on the customer's preferences, which are unknown to the leader. We next analyze how these responses can be determined.

\subsection{The Followers' Decision Problem}\label{sect:follower}
{For a given leader solution  ${\bar{\bm{\ell}}}$, 
customers who are not offered home delivery face a \textit{relatively simple} decision: identifying their preferred collection point among those activated by the leader. In contrast, customers who are offered home delivery face a \textit{more complex} decision. Specifically, they must decide whether to accept the offered delivery time. If the offer is rejected, they must then identify their preferred collection point among those activated by the leader.}\\

In order to decide whether to accept the leader's offer, such customers must compare the proposed service times with their own preference intervals. For this, for each customer $i \in I$, we define the following decision variables:

\begin{itemize}
\item[] $z_{ir} =
\begin{cases}
1 & \text{if customer $i \in I$ is offered home delivery in the time interval }\Delta_{i\,r}, r\in R^i,\\
0 & \text{otherwise.}
\end{cases}
$
\end{itemize}

\begin{itemize}
\item[] $\beta_{ik} =
\begin{cases}
1 & \text{if customer $i \in I$ picks up the parcel at collection point $k \in K$}, \\
0 & \text{otherwise.}
\end{cases}$

\item[] $\alpha_{ir} =
\begin{cases}
1 & \text{if customer $i \in I$ accepts the service time offered by the leader,} \\
& \text{which belongs to their time interval $r \in R^i$}, \\
0 & \text{otherwise.}
\end{cases}$
\end{itemize}

Then, the follower's problem for a given customer $i \in I$ can be modeled as follows:

\begin{subequations}
\begin{align}
F(\bm \ell):\quad  \min_{\bm \alpha, \bm \beta, \mathbf{z}} &  \sum_{i \in I} \sum_{r \in R^i} {\gamma}_{ir} \alpha_{ir} + \sum_{i \in I} \sum_{k \in K} {\gamma}_{i,k+|R^i|} \beta_{ik}    \label{eq:F_fobj-0}&&\\
\text{subject to:\quad} & \sum_{r\in R^i} \alpha_{ir}+\sum_{k\in K}\beta_{ik} = 1 && i\in I \label{eq:F_ir_a_collection-0}\\
  & \beta_{ik}\leq \bar{y}_{k} \qquad && i\in I, k\in K \label{eq:F_recoger_abierta-0}\\
& \alpha_{ir}\leq {z}_{ir} \qquad && i\in I, r\in R^i \label{eq:F_recibir_casa-0}\\
&  \sum_{r\in R^i} z_{ir}= \bar s_i \qquad\qquad && i\in I\label{eq:served-0-sub} \\
&  {\sum_{r\in R^i} \left(\delta_{i,\, r-1} +\varepsilon\right) z_{ir}\leq \bar g_i}\qquad\qquad && i\in I\label{eq:interval1-0} \\
&   {\bar g_i \leq \sum_{r\in R^i} \delta_{i, r}\, z_{ir}}\qquad\qquad && i\in I \label{eq:interval2-0} \\
 & \alpha_{ir} \in\{0, 1\}\qquad && i\in I, r\in R^i \label{eq:F_alpha-0} \\
 & \beta_{ik}\in\{0, 1\}\qquad && i\in I, k\in K \label{eq:F_beta-0}\\
  &   z_{ir}\in \{0,1\}  \qquad \qquad && i\in I, r\in R^i. \label{eq:L_zvar-0}
\end{align}
\end{subequations}

The objective function \eqref{eq:F_fobj-0} minimizes customers’ preferences for receiving or picking up their parcels. Constraints \eqref{eq:F_ir_a_collection-0} ensure that each customer selects a single service mode: either the proposed home delivery time or a collection point for parcel pickup. Constraint \eqref{eq:F_recoger_abierta-0} prevents any customer from picking up their parcel at a collection point that has not been activated at the upper level. Clearly, customers may accept the proposed home delivery only if they are included in the service route, see \eqref{eq:F_recibir_casa-0}-\eqref{eq:served-0-sub}. Note that \eqref{eq:F_recibir_casa-0} ensures that if a customer $i$ accepts home delivery in an interval with preference $\gamma_{ir}$, then the service time offered by the leader lies within the interval {$\Delta_{ir}$, i.e., $z_{ir}=1$ for some index $r\in R^i$.}
The constraints \eqref{eq:interval1-0}-\eqref{eq:interval2-0} identify the 
interval containing the service time proposed by the leader ($\bar g_i$), where $\varepsilon$ is a small positive value representing a left-open interval. Finally, constraints \eqref{eq:F_alpha-0}-\eqref{eq:L_zvar-0} enforce the binary nature of the customer decision variables. 

\begin{rem}
 We observe that the group of constraints \eqref{eq:served-0-sub}-\eqref{eq:interval2-0}, together with the dependent variables {$z_{ir}$, $i\in I$, $r\in R^i$,} 
 can be transferred to the leader's problem. This has the clear advantage of simplifying the formulation of the follower’s problem, which, as we will see next, can be solved and dualized using classical linear programming techniques.\\
 In the following we denote by $\mathcal{X}^{\bm z}(\mathbf s, \mathbf g)$ the domain determined by binary vectors $\bf z$ that satisfy constraints \eqref{eq:served-0-sub}-\eqref{eq:interval2-0} for a given vector $(\mathbf s,\, \mathbf g)$.\\
\end{rem}

When the variables {$\bar z_{ir}$, $i \in I$, $r \in R^i$}, determined by a leader's solution ${\bar{\bm{\ell}}}$ are also obtained from the leader's problem, then the followers' problem corresponding to a solution $({\bar{\bm{\ell}}}, {\bar{\mathbf{z}}})$ reduces to:
\begin{subequations}
\begin{align}
F({\bar{\bm{\ell}}}, {\bar{\mathbf{z}}}):\quad  \min_{\bm \alpha, \bm \beta} &  \sum_{i \in I} \sum_{r \in R^i} {\gamma}_{ir} \alpha_{ir} + \sum_{i \in I} \sum_{k \in K} {\gamma}_{i,k+|R^i|} \beta_{ik}    \label{eq:F_fobj-00}&&\\
\text{subject to:\quad} & \sum_{r\in R^i} \alpha_{ir}+\sum_{k\in K}\beta_{ik} = 1 && i\in I \label{eq:F_ir_a_collection-00}\\
  & \beta_{ik}\leq \bar{y}_{k} \qquad && i\in I, k\in K \label{eq:F_recoger_abierta-00}\\
& {\alpha_{ir}\leq {\bar{z}_{ir}}} \qquad && i\in I, r\in R^i \label{eq:F_recibir_casa-00}\\
 & \alpha_{ir} \in\{0, 1\}\qquad && i\in I, r\in R^i \label{eq:F_alpha-00} \\
 & \beta_{ik}\in\{0, 1\}\qquad && i\in I, k\in K. \label{eq:F_beta-00}
\end{align}
\end{subequations}

\begin{rem}{$F({\bar{\bm{\ell}}}, {\bar{\mathbf{z}}})$} can be decomposed in $|I|$ independent problems, one for each customer. 
Since all coefficients are binary (0 or 1) and the right-hand side equals 1, the integrality restrictions on the variables can be relaxed. Hence, the follower’s problem for each customer $i \in I$ can be written as follows:
\begin{subequations}
\begin{align}
F^i({\bar{\bm{y}}}, {\bar{\mathbf{z}}}):\quad  \min_{\bm \alpha ^i, \bm \beta ^i} &  \sum_{r \in R^i} {\gamma}_{ir} \alpha_{ir} + \sum_{k \in K} {\gamma}_{i,k+|R^i|} \beta_{ik}    \label{eq:F_fobj-i-sub}&&\\
\text{subject to:\quad} & (\bm\alpha^i, \bm \beta^i)\in\mathcal X^{\mathbf\alpha^i, \bm\beta^i}(\bm{\bar y}, \mathbf{\bar z}),\nonumber
\end{align}
\end{subequations} 

where $\mathcal X^{\mathbf\alpha^i, \bm\beta^i}(\bm{\bar y}, \mathbf{\bar z})=\{(\bm\alpha^i, \bm \beta^i) \text{ that satisfy } \eqref{eq:F_ir_a_collection-00}–\eqref{eq:F_recibir_casa-00},\, \alpha_{ir}\geq 0, \, r\in R^i; \beta_{ik}\geq 0,\, k\in K\}$.\\

Note that the values of the vector $\bm \beta$, allow the leader to reproduce the followers’ responses, not only for customers served at activated collection  points but also for customers who accept home delivery at the proposed service time. In particular, for all $i \in I$, the value $\sum_{k \in K} \beta_{ik}$ serves as an indicator of the response of customer $i \in I$: it takes the value 0 if customer $i \in I$ is visited by the leader’s tour and the proposed service time is accepted, and 1 otherwise. Hence, the compensation received by a customer who accepts home delivery at the proposed service time can be expressed as $\theta^T_i h_i \left(1 - \sum_{k \in K} \beta_{ik}\right)$.\\
\end{rem}

\subsection{Modeling Cost Compensations}

As mentioned above, the time-related compensation paid to customers who accept home delivery at the proposed service time depends on the deviation between the actual and ideal service times. For a given customer $i\in I$, it can be expressed as $\theta_i^T h_i \left(1-\sum_{k\in K}\beta_{ik}\right)$.\\

Determining the compensation for a customer who is not offered home delivery or who chooses to pick up their parcel at a collection point is more complex, as it depends on the distance to the closest activated collection point, which may differ from the collection point ultimately selected by the customer. 
A closed-form of this compensation for a given customer $i\in I$,  is $\left(\sum_{k\in K}\beta_{ik}\right)\left(\min_{k\in K}\{\theta^F_{ik}y_k\}\right)$.  Next, we introduce an additional set of decision variables that allows us to derive a linear system for this expression.\\

For all $i\in I$, $k\in K$, let \\ $$\rho_{ik}=\begin{cases}
1 & \text{if $k\in K$ is the activated collection point closest to customer } i \in I,\\ & \text{who does not receive home delivery},\\
0 & \text{otherwise.}
\end{cases}
$$
The following set of constraints determines the correct values of variables $\rho_{ik}$, $i\in I$, $k\in K$:  
\begin{subequations}
\begin{align}
& \sum_{k\in K}\rho_{ik}=1-\sum_{r\in R^i}\alpha_{ir}\qquad && i\in I,\label{eq:i_not-home}\\
& \sum_{\substack{k':\\d_{ik'}>d_{ik}}}\rho_{ik'}+y_{k}\leq 1&& i\in I,\ k\in K\label{eq_closer than-k}\\
& \rho_{ik}\leq y_k&& i\in I,\, k\in K. \label{eq_open-facility}
\end{align}
\end{subequations}

For all $i\in I$, the right-hand side of \eqref{eq:i_not-home} takes the value 1 only if customer $i$ does not receive home delivery. Thus, this constraint activates one variable $\rho_{ik}$ for such customers. By \eqref{eq_closer than-k}, if collection point $k \in K$ is open, then, for any customer $i\in I$, no variable $\rho_{ik'}$ associated with a collection point farther from customer $i$ than collection point $k$ can be activated.
Finally, \eqref{eq_open-facility} ensure that  $\rho_{ik}=1$ only if collection point $k$ is open.\\

Indeed, the compensation paid to customers served through a collection point can now be expressed as:

$$\sum_{i\in I}\sum_{k\in K}\theta_{ik}^F\rho_{ik}.$$

Let $\mathcal X^{\bm\rho}(\bm{\bar y}, \bm{ \bar \alpha})=\{\bm\rho: \text{ that satisfy } \eqref{eq:i_not-home}-\eqref{eq_open-facility}\}$ for a given vector $\left(\bm{\bar y}, \bm{\bar \alpha}\right)$.

\subsection{A Mixed-Integer Bilevel Formulation}

The considerations described above lead to the following mixed-integer bilevel program (MIBP) for the \ABC: 

\begin{subequations} \begin{align} \min_{\mathbf y,\mathbf s,\mathbf x, \mathbf g, \mathbf h, \mathbf z,\bm \alpha, \bm \beta, \bm \rho}\quad & \sum_{k\in K} f_k y_k + \sum_{(i,j)\in A} c_{ij}x_{ij} + \sum_{i\in I}\theta_i^T h_i \left(1-\sum_{k\in K}\beta_{ik}\right) + \sum_{i\in I}\sum_{k\in K}\theta_{ik}^F\rho_{ik} \label{eq:bilevel_obj}\\ 
\text{subject to}\quad & (\mathbf y, \mathbf s, \mathbf x, \mathbf g, \mathbf h)\in \mathcal X^{\bm\ell} \nonumber \\ 
& \mathbf z\in \mathcal X^{\bm z}(\mathbf s, \mathbf g) \nonumber \\ 
& (\bm\alpha^i,\bm\beta^i) \in \underset{\left(\bm \alpha^i,\,\bm \beta^i\right)\in \mathcal X^{\mathbf\alpha^i, \bm\beta^i}(\mathbf{\bar y}, \mathbf{\bar z})}{\arg\min} \left\{ \sum_{r\in R^i}\gamma_{ir}\alpha_{ir} + \sum_{k\in K}\gamma_{i,k+|R^i|}\beta_{ik} \right\}, && i\in I, \label{eq:bilevel_follower} \\
&\bm\rho\in \mathcal X^{\bm\rho}(\bm{y}, \bm{\alpha}).&& \label{comput_rho}
\end{align} \end{subequations}

The leader's objective function \eqref{eq:bilevel_obj} minimizes the total system cost, consisting of four components: the cost of opening collection points, the routing cost, the compensation paid to customers who accept home delivery at the proposed service time, and the compensation associated with customers who pick up their parcels at collection points.
As explained above, the term $\theta_i^T h_i \left(1-\sum_{k\in K}\beta_{ik}\right)$ accounts for the time-related compensation paid to customer $i\in I$ only when home delivery is accepted, whereas $\theta_{i\bar k}^F\rho_{i\bar k}=\sum_{k\in K}\theta_{ik}^F\rho_{ik}$ represents the compensation paid when customer $i$ picks up their parcel at some collection point, where $\bar k \in K$ denotes the activated collection point closest to customer $i$. Imposing that the decision variables belong to the domains determined by $\mathcal X^{\bm\ell}$ and $\mathcal X^{\bm z}(\mathbf s, \mathbf g)$ guarantee the nature and domains of the leader's decision variables. Additionally, it is worth highlighting that constraint \eqref{eq:bilevel_follower} ensures that, for each customer $i\in I$, the response $(\bm{\alpha}^i,\bm{\beta}^i)$ corresponds to an optimal solution of problem 
$F^i({\bar{\bm{y}}}, {\bar{\mathbf{z}}})$ parameterized by the leader's decisions and \eqref{comput_rho} ensures that the collection points closest to the customers who do not receive home delivery are identified  correctly.\\

For this bilevel model, the classical optimistic approach is assumed \citep{kalashnikov2016bilevel, dempe2020bilevel}. This assumption is reflected in \eqref{eq:bilevel_obj}, where the follower's decision variables are included in the leader's objective function. Consequently, when multiple optimal solutions to the follower's problem exist for a given upper-level decision, the response that is most favorable to the company (leader) is selected. This assumption is appropriate in the context of the problem because customers and the company do not have inherently conflicting objectives. Customers select their most preferred service alternative among those made available by the company, and, in the case of indifference, the company may guide the final service choice without worsening the customer's preference level. Therefore, the optimistic convention does not require a customer to accept a less-preferred alternative; it only resolves ties among service alternatives that are equally preferred by that customer.

\section{An Equivalent Single-Level Formulation}\label{sect:single-level} 
In order to develop a single-level reformulation, we apply the primal-dual optimality conditions to the followers' problem.  \\

By noting that, for a customer $i\in I$, constraints \eqref{eq:F_ir_a_collection-00}-\eqref{eq:F_recibir_casa-00} can be expressed as a single constraint, namely $\sum_{r\in R^i} \bar z_{ir}\alpha_{ir}+\sum_{k\in K}\bar y_k\beta_{ik} = 1$,
the dual of $F^i(\bm{\bar y}, \bm{\bar{z}})$ 
for customer $i \in I$ is:
\begin{subequations}
\begin{align}
D^i(\bm{\bar y}, \bm{\bar z}):\quad  \max_{\bm u} &\quad   u_i \label{eq:F_fobj-i-dual}&&\\
\text{subject to:\quad} & \bar z_{ir} u_i \leq \gamma_{ir} && r\in R^i, \label{eq:dual1}\\
&\bar y_{k} u_i \leq \gamma_{i, k+|R^i|} && k\in K. \label{eq:dual2}
\end{align}
\end{subequations}

Thus, a classical single-level reformulation of the MIBP proposed for the \ABC is given by:

\begin{subequations}
\begin{align}
\text{SLR}:\qquad \min_{\bm y,\bm s,\bm x,\bm g,\bm h,\bm z,\bm{\alpha}, \bm{\beta}, \bm\rho, \bm u} \qquad  &  \sum_{k\in K}f_ky_k + \sum_{(i,j)\in A} c_{ij} x_{ij}+ \sum_{i \in I} \theta^T_i h_i&&\hspace{-10.pt}\left(1-\sum_{k\in K}\beta_{ik}\right)\,\nonumber\\
& &&+ \sum_{i \in I}\sum_{k \in K} \theta^F_{ik}\rho_{ik}  \tag{\ref*{eq:bilevel_obj}}\\
\text{subject to}\quad 
 & \eqref{eq:dual1}-\eqref{eq:dual2}\qquad && i\in I,\nonumber\\
& u_i\geq \sum_{r\in R^i}\gamma_{ir}\alpha_{ir}+\sum_{k\in K} {\gamma}_{i,k+|R^i|} \beta_{ik} 
&& i\in I,   \label{pd-1}\\
& (\mathbf y, \mathbf s, \mathbf x, \mathbf g, \mathbf h)\in \mathcal X^{\bm\ell}, \nonumber \\ 
& \mathbf z\in \mathcal X^{\bm z}(\mathbf s,\,\mathbf g),\nonumber \\
& \left(\bm \alpha^i,\,\bm \beta^i\right)\in \mathcal X^{\mathbf\alpha^i, \bm\beta^i}(\mathbf{\bar y}, \mathbf{\bar z}) && i\in I,\nonumber\\
&\bm\rho\in \mathcal X^{\bm\rho}(\bm{y}, \bm{\alpha}).&&\nonumber
\end{align}
\end{subequations}

\noindent where constraints \eqref{eq:dual1}-\eqref{eq:dual2}
together with \eqref{pd-1}, and $\left(\bm \alpha^i,\,\bm \beta^i\right)\in \mathcal X^{\mathbf\alpha^i, \bm\beta^i}(\mathbf{\bar y}, \mathbf{\bar z})$, for all $i\in I$, model the primal-dual optimality conditions for the followers' problems and thus ensure that $\bm{\alpha}$ and $\bm{\beta}$ correspond to their optimal responses.

\subsection{Linearization of the SLR} \label{sec:MILP-I}
As can be seen, the SLR involves nonlinear terms both in the objective function and in some of the constraints, e.g., 
\eqref{eq:dual1}-\eqref{eq:dual2}. 
This prevents the formulation from being handled efficiently from a computational perspective. Next, we describe how these nonlinearities can be linearized.\\

\begin{itemize}
\item  In the objective function, the term $\sum_{i \in I} \theta^T_i h_i \left(1-\sum_{k \in K}\beta_{ik}\right)$ can be linearized by introducing an auxiliary set of decision variables $\bar h_i$ to represent the actual deviation from the ideal service time for customers who accept home delivery. The values of these variables depend on the followers' responses: for customers $i \in I$ who accept the leader's offer, $\bar h_i$ takes the value $h_i$, and 0 otherwise. In particular, for each customer, let:

\begin{itemize}
\item[]$\bar h_i=\begin{cases}h_i &\text{if customer $i\in I$ is offered home delivery and accepts it,}\\
0 & \text{otherwise}.\end{cases}$
\end{itemize}

\noindent The values of these auxiliary variables can be determined by the following inequalities:
  \begin{subequations}
  \begin{align}
  & \bar h_i\leq h_i \qquad && i\in I,\label{h_bar_1}\\
   & \bar h_i\leq  \overline{\text{M}}_i \left(1-\sum_{k\in K}\beta_{ik}\right) \qquad&& i\in I,\label{h_bar_2}\\
  & \bar h_i\geq h_i- \overline{\text{M}}_i\sum_{k\in K}\beta_{ik} \qquad && i\in I,\label{h_bar_3} 
  \end{align}
  \end{subequations}
  \noindent where, for all $i\in I$, $\overline{\text{M}}_i=\max\{\delta_{i\bar r^i}-T^*_i, T^*_i-\delta_{i0}\}$ is an upper bound on the maximum deviation of its actual service time from its ideal service time.
\item Constraints \eqref{eq:dual1}-\eqref{eq:dual2} 
can be linearized using standard techniques as:
\begin{subequations} 
\begin{align} 
u_i 
&\leq \gamma_{ir} 
+ H_{ir}\left(1-z_{ir}\right)
&& i \in I,\ r \in R^i,
\label{eq:dual_2_linearized}
\\[0.5em]
u_i 
&\leq \gamma_{i,k+|R^i|} 
+ H_{i, k+|R^i|}(1-y_k)
&& i \in I,\ k \in K,
\label{eq:dual_1_linearized}
\end{align}
\end{subequations}
\noindent where, in each case, $H_{il}=|L^i|-\gamma_{il}$, $l\in L^i$ is a suitable upper bound on the maximum value of $u_i$.

\end{itemize}

Therefore, replacing the nonlinear terms in the SLR with the linearization constraints described above yields the following mixed-integer linear reformulation, hereafter referred to as MILP-I:

\begin{subequations} \begin{align} \text{MILP-I}:\quad \min_{\bm y,\bm s,\bm x,\bm g,\bm h,\bm{\bar h}, \bm z,\bm{\alpha},\bm{\beta}, \bm \rho,\bm u}\quad & \sum_{k\in K}f_k y_k + \sum_{(i,j)\in A}c_{ij}x_{ij} + \sum_{i\in I}\theta_i^T\bar h_i + && \sum_{i\in I}\sum_{k\in K}\theta_{ik}^F\rho_{ik} \label{eq:MILP-I-obj}\\
\text{subject to}\quad 
& \eqref{h_bar_1}-\eqref{h_bar_3} && i\in I,  \nonumber\\ 
& \eqref{eq:dual_2_linearized}-\eqref{eq:dual_1_linearized} && i\in I,  \nonumber\\  
&\bar h_i\geq 0 \qquad \qquad &&   i\in I,  \nonumber\\ 
& (\mathbf y, \mathbf s, \mathbf x, \mathbf g, \mathbf h)\in \mathcal X^{\bm\ell}, \nonumber \\ 
& \mathbf z\in \mathcal X^{\bm z}(\mathbf s,\,\mathbf g) \nonumber,\\ 
& \left(\bm \alpha^i,\,\bm \beta^i\right)\in \mathcal X^{\mathbf\alpha^i, \bm\beta^i}(\mathbf{\bar y}, \mathbf{\bar z}) && i\in I,\nonumber\\
&\bm\rho\in \mathcal X^{\bm\rho}(\bm{y}, \bm{\alpha}).&&\nonumber\end{align} 
\end{subequations}

\subsection{An Alternative Single-Level Reformulation Using Only Primal Variables} \label{sec:MILP-II}

Below, we present an alternative linearization to the one described above, in which the primal-dual optimality constraints \eqref{eq:dual1}-\eqref{eq:dual2} and \eqref{pd-1} are replaced with equivalent linear inequalities involving only primal variables. To this end, we introduce the following notation.\\

Let $K^z_{\prec}(i,r)=\{k\in K: \gamma_{i, k+|R^i|}<\gamma_{ir}\}$ be the set of collection points that would be preferred by customer $i\in I$ to home delivery with service time in the interval $\Delta_{ir}$.\\

Let also $K^y_{\succ}(i, k)=\{k'\in K:  \gamma_{i,k'+|R^i|}>\gamma_{i,k+|R^i|} \}$ be the set of collection points less preferred than collection point $k \in K$ for customer $i \in I$.\\

Consider the following system of linear inequalities:
\begin{subequations}
\begin{align}
& {z_{ir}}  \leq \alpha_{ir}+\sum_{k\in K^z_{\prec}(i,r)}\beta_{ik}\qquad && i\in I, r\in R^i, \label{SL2-0} \\
&\alpha_{ir} +y_k\leq 1\qquad && i\in I, r\in R^i, k\in K^z_{\prec}(i,r), \label{SL2-1} \\
& y_k+\sum_{k'\in K^y_{\succ}(i,k)}\beta_{ik'}\leq 1 && i\in I, k\in K, \label{SL2-2}\\
&(\bm\alpha^i, \bm\beta^i)\in\mathcal X^{\bm\alpha^i, \bm\beta^i}(\bm{ y}, \mathbf{z})&& i\in I.\nonumber
\end{align}
\end{subequations}

\begin{proposition}\label{propo:primal}
Let $(\bm{\bar\ell}, \mathbf{\bar z})$, with $\bm{\bar \ell}=(\mathbf{\bar y},  \mathbf{\bar s}, \mathbf{\bar x}, \mathbf{\bar g}, \mathbf{\bar h})\in\mathcal X^{\bm\ell}$ and $\mathbf{\bar z}\in\mathcal X^{\bm z}(\mathbf{\bar s}, \mathbf{\bar g})$, be a given upper-level solution. Then, any solution $(\bm{\alpha}, \bm{\beta})$ with $(\bm{\alpha}^i, \bm{\beta}^i)\in \mathcal X^{\bm\alpha^i, \bm\beta^i}(\mathbf{\bar y}, \mathbf{\bar z})$ for all $i\in I$, that satisfies 
\eqref{SL2-0}-\eqref{SL2-2} determines an optimal solution to subproblem $F(\bm{\bar\ell}, \mathbf{\bar z})$.
\end{proposition}
\proof
Let $(\bm{\bar\alpha^i}, \bm{\bar\beta^i})\in \mathcal X^{\bm\alpha^i, \bm\beta^i}(\mathbf{\bar y}, \mathbf{\bar z})$ denote a solution to the subsystem defined by 
\eqref{SL2-0}-\eqref{SL2-2}, associated with $(\bm{\bar y}, \mathbf{\bar z})$ for a given $i\in I$. We now show that $(\bm{\bar\alpha^i}, \bm{\bar\beta^i})$ is an optimal solution to $F^i(\bm{\bar y}, \mathbf{\bar z})$.

We proceed by contradiction.

Suppose, on the contrary, that $(\bm{\bar\alpha^i}, \bm{\bar\beta^i})$ is not optimal for $F^i(\bm{\bar y}, \mathbf{\bar z})$.
Then, there exists a solution to the subsystem associated with $i$, $(\bm{{\alpha^i}^*}, \bm{{\beta^i}^*})\in \mathcal X^{\bm\alpha^i, \bm\beta^i}(\mathbf{\bar y}, \mathbf{\bar z})$, such that,
$$\sum_{r \in R^i} {\gamma}_{ir} \alpha_{ir}^* + \sum_{k \in K} {\gamma}_{i,k+|R^i|} \beta_{ik}^*<\sum_{r \in R^i} {\gamma}_{ir} \bar\alpha_{ir} + \sum_{k \in K} {\gamma}_{i,k+|R^i|} \bar\beta_{ik}.$$

Since both $(\bm{{\bar\alpha^i}}, \bm{{\bar \beta^i}}),\,(\bm{{\alpha^i}^*}, \bm{{ \beta^i}^*})\in \mathcal X^{\bm\alpha^i, \bm\beta^i}(\mathbf{\bar y}, \mathbf{\bar z})$, 
in each case exactly one component takes the value 1, while all remaining components are equal to 0. In particular, exactly one of the following cases occurs:
\begin{itemize}
\item There exist $\bar r, r^*\in R^i$ such that $\bar\alpha_{i\bar r}=\alpha^*_{ir^*}=1$, with $\bar z_{i\bar r}=\bar z_{ir^*}=1$ and $\gamma_{ir^*}<\gamma_{i\bar r}$.\\ 
This would imply that $\bar r\ne r^*$, and hence $\sum_{r\in R^i}\bar z_{ir}= 2$. Therefore, $(\bm{{\bar\alpha^i}}, \bm{{\bar \beta^i}})$ would violate constraint \eqref{eq:served-0-sub} corresponding to $i$. Hence, this case is not possible. 
\item There exist $\bar r\in R^i$ and $k^*\in K$ such that $\bar\alpha_{i\bar r}=1$ and $\beta^*_{ik^*}=1$, with $\bar z_{i\bar r}=\bar y_{k^*}=1$ and $\gamma_{i, k^*+|R^i|}<\gamma_{i\bar r}$.\\
Since $\gamma_{i, k^*+|R^i|}<\gamma_{i\bar r}$, it follows that $k^*\in K^z_{\prec}(i, \bar r)$, with $\bar\alpha_{i\bar r}=\bar y_{k^*}=1$. Therefore, this case is not possible, since constraint \eqref{SL2-1} would be violated by $(\bm{{\bar\alpha^i}}, \bm{{\bar \beta^i}})$.

\item There exist $\bar k\in K$ and $r^*\in R^i$ such that $\bar\beta_{i\bar k}=1$ and $\alpha^*_{i r^*}=1$, with $\bar y_{\bar k}=\bar z_{ir^*}=1$ and $\gamma_{ir^*}<\gamma_{i, \bar k+|R^i|}$. \\
Since $\gamma_{ir^*}<\gamma_{i, \bar k+|R^i|}$, it follows that $\bar k\notin K_{\succ}^z(i, r^*)$, and hence $\sum_{k\in K_{\succ}^z(i, r^*)} \bar \beta_{ik}=0$. Therefore, $(\bm{{\bar\alpha^i}}, \bm{{\bar \beta^i}})$ would violate the constraint \eqref{SL2-0} for index $r^*$, since $\bar z_{ir^*}=1$ and $\bar \alpha_{i r^*}=0$.

\item There exist $\bar k, k^*\in K$ such that $\bar\beta_{i\bar k}=\beta^*_{i\bar k^*}=1$, with $\bar y_{\bar k}=\bar y_{k^*}=1$ and $\gamma_{i, k^*+|R^i|}<\gamma_{i, \bar k+|R^i|}$.\\
Since $\gamma_{i, k^*+|R^i|}<\gamma_{i, \bar k+|R^i|}$, it follows that $\bar k\in K^y_{\succ}(i, k^*)$, with $\bar y_{k^*}=\bar\beta_{i\bar k}=1$. Therefore, this case is not possible, since $(\bm{{\bar\alpha^i}}, \bm{{\bar \beta^i}})$ would violate constraint \eqref{SL2-2}.
\end{itemize}

Therefore, using constraints \eqref{SL2-0}-\eqref{SL2-2} to characterize the followers’ optimal responses, while maintaining the linearization \eqref{h_bar_1}-\eqref{h_bar_3} of the upper-level objective function, we obtain the following alternative mixed-integer linear reformulation, referred to as MILP-II:

\begin{subequations} \begin{align} \text{MILP-II}:\quad \min_{\bm y,\bm s,\bm x,\bm g,\bm h,\bm{\bar h},\bm z,\bm{\alpha},\bm{\beta},\bm\rho}\quad & \sum_{k\in K}f_k y_k + \sum_{(i,j)\in A}c_{ij}x_{ij} + \sum_{i\in I}\theta_i^T\bar h_i + && \sum_{i\in I}\sum_{k\in K}\theta_{ik}^F\rho_{ik} \label{eq:MILP-I-obj}\\ 
\text{subject to}\quad 
& \eqref{h_bar_1}-\eqref{h_bar_3} && i\in I,  \nonumber\\ 
&\bar h_i\geq 0 \qquad \qquad &&   i\in I,  \nonumber\\ 
& \eqref{SL2-0}-\eqref{SL2-2}, \nonumber\\
& (\mathbf y, \mathbf s, \mathbf x, \mathbf g, \mathbf h)\in \mathcal X^{\bm\ell}, \nonumber \\ 
& \mathbf z\in \mathcal X^{\bm z}(\mathbf s,\,\mathbf g), \nonumber \\
& \left(\bm \alpha^i,\,\bm \beta^i\right)\in \mathcal X^{\mathbf\alpha^i, \bm\beta^i}(\bm{y}, \mathbf{z}) && i\in I,\nonumber\\
&\bm\rho\in \mathcal X^{\bm\rho}(\bm{y}, \bm{\alpha}).&&\nonumber
\end{align} 
\end{subequations}

\section{A Mixed-Integer Programming Heuristic for the \ABC} \label{sec:heur}

The proposed Mixed-Integer Programming (MIP) heuristic follows a two-phase structure designed to construct high-quality feasible solutions that can be used as warm starts for the single-level MILP formulations proposed in Sections \ref{sec:MILP-I} and \ref{sec:MILP-II}. Phase 1 solves an approximation, consisting of a reduced location-routing problem considering only the depot and the candidate collection points. The objective of this phase is to identify a promising subset of collection points, together with initial customer-to-collection-point assignments and a routing structure connecting the selected collection points. In Phase 2, this partial initial solution is used as input to a parameterized, leader-oriented version of the MILP formulation, in which only the collection points selected in Phase 1 may be opened. In addition, customer rejection is not allowed; thus, each customer must be served through home delivery, if proposed by the leader. This auxiliary formulation reduces the search space and constructs a complete feasible solution that can subsequently be used as an incumbent for the original MILP formulation.\\

\textbf{Phase 1. Reduced Location-Routing Problem.} The first phase solves a simplified location-routing problem that considers only the candidate collection points and the depot. The purpose of this phase is to identify a promising subset of collection points and construct a proxy tour connecting them. Home-delivery decisions, time-window decisions, and customers' time-related preferences are not considered at this stage. However, the model accounts for customers' preferences over collection points.\\

Let $D=\{0\}$ represent the depot, and define $N^{L}=K\cup D$ as the node set considered in Phase 1. The corresponding arc set is $A^{L}=\{(i,j): i,j\in N^{L}, i\neq j\}$. The decision variables are defined as follows:

\begin{itemize}[label={}]
\item $ y_k^{L} =
\begin{cases}
1, & \text{if collection point } k \text{ is opened in Phase 1},\\
0, & \text{otherwise},
\end{cases}
\qquad k\in K,$

\item $\beta_{ik}^{L} =
\begin{cases}
1, & \text{if customer } i \text{ is assigned to collection point } k,\\
0, & \text{otherwise},
\end{cases}
\qquad i\in I,\ k\in K, $

\item $ x_{ij}^{L} =
\begin{cases}
1, & \text{if arc } (i,j) \text{ is used in the proxy tour},\\
0, & \text{otherwise},
\end{cases}
\qquad (i,j)\in A^{L}, $

\item $\rho_{ik}^{L} =
\begin{cases}
1, & \text{if the cost compensation of customer } i \text{ is } \theta_{ik}^F,\\
0, & \text{otherwise},
\end{cases}
\qquad i\in I,\ k\in K, $

\item and $v_k^{L}$ is an auxiliary variable associated with collection point $k\in K$, used in the MTZ subtour-elimination constraints.

\end{itemize}

The Phase 1 model is formulated as follows:
\begin{subequations}
\begin{align}
\min_{\mathbf{y},\bm{\beta},\bm{\rho},\mathbf{x}, \mathbf{v}} \quad
& \sum_{k\in K} f_k y_k^{L}
+ \sum_{(i,j)\in A^{L}} c_{ij}x_{ij}^{L}+ \sum_{i\in I}\sum_{k\in K} \theta_{ik}^{F}\rho_{ik}^{L}
\label{eq:phase1_obj}
\\[0.2cm]
\text{s.t.}\quad
& \beta_{ik}^{L} \leq y_k^{L}
&& i\in I,\; k\in K,
\label{eq:phase1_assign_open}
\\
& \sum_{k\in K} \beta_{ik}^{L} = 1
&& i\in I,
\label{eq:phase1_assign_customer}
\\
& \sum_{\substack{j\in N^{L}\\j\neq 0}} x_{0j}^{L} = 1
\label{eq:phase1_leave_depot}
\\
& \sum_{\substack{j\in N^{L}\\j\neq 0}} x_{j0}^{L} = 1
\label{eq:phase1_return_depot}
\\
& \sum_{\substack{j\in N^{L}\\j\neq k}} x_{kj}^{L} = y_k^{L}
&& k\in K,
\label{eq:phase1_flow_out}
\\
& \sum_{\substack{j\in N^{L}\\j\neq k}} x_{jk}^{L} = y_k^{L}
&& k\in K,
\label{eq:phase1_flow_in}
\\
& v_i^{L} - v_j^{L} + |K|x_{ij}^{L} \leq |K|-1
&& i,j\in K,\; i\neq j,
\label{eq:phase1_mtz}
\\
& y_k^{L}
+ \sum_{k'\in K^{y}(i,k)} \beta_{ik'}^{L}
\leq 1
&& i\in I,\; k\in K,
\label{eq:phase1_prop1}
\\
& \sum_{k\in K}\rho_{ik}^L=1\qquad && i\in I,\label{eq:i_not-home-Ph1}\\
& \sum_{\substack{k':\\d_{ik'}>d_{ik}}}\rho_{ik'}^L+y_{k}^L\leq 1 && i\in I,\ k\in K,\label{eq_closer than-k-Ph-1}\\
& \rho_{ik}^L\leq y_k^L && i\in I,\, k\in K, \label{eq_open-facility-Ph1}\\
& y_k^{L}\in\{0,1\}
&& k\in K,
\\
& \beta_{ik}^{L}\in\{0,1\}
&& i\in I,\; k\in K,
\\
& \rho_{ik}^{L}\in\{0,1\}
&& i\in I,\; k\in K,
\\
& x_{ij}^{L}\in\{0,1\}
&& (i,j)\in A^{L},
\\
& v_k^{L} \geq 0
&& k\in K.
\end{align}    
\end{subequations}

 The solution to the Phase 1 problem identifies a collection-point configuration and a customer assignment that is optimal in terms of customer preferences for the current configuration, as imposed by constraints \eqref{eq:phase1_prop1}, while also accounting for collection point-activation costs, routing costs, and the compensation costs associated with the assignment of customers to closest activated collection points. Note that now the right-hand-side of \eqref{eq:i_not-home-Ph1} is 1, given that all customers are served at a collection point.\\
 
 Let $(\mathbf{\hat y}^L,\bm{\hat\beta}^L, \bm{\hat\rho}^L, \mathbf{\hat x}^L,\mathbf{\hat v}^L)$ denote an optimal solution to the Phase 1 problem. 
This solution specifies the collection points that are allowed to open in Phase 2, namely those such that $\hat y_{k}^L=1$ and provides an initial solution for Phase 2. Specifically, the values of $\mathbf{\hat y}^L$, $\mathbf{\hat x}^L$, and $\bm{\hat\beta}^L$ initialize the corresponding opening, routing, and assignment variables, respectively. Overall, these decisions represent an initial service configuration in which customers are served exclusively through collection points, with the route visiting only the activated facilities.\\

\textbf{Phase 2: Restricted routing formulation}

Using 
the partial solution $(\mathbf{\hat y}^L,\bm{\hat\beta}^L, \bm{\hat\rho}^L, \mathbf{\hat x}^L,\mathbf{\hat v}^L)$ obtained in Phase 1, Phase 2 solves an auxiliary  formulation restricted to the activated collection points, in which customers 
can be routed although
the possibility of rejecting the service option proposed by the leader is temporarily disregarded. Let $\mathcal{X}^{\mathrm{MILP-II}}$ denote the feasible region induced by the MILP-II formulation. The resulting Phase 2 model is formulated as follows:

\begin{subequations}
\begin{align}
\qquad
\min \quad
& \Phi^{\mathrm{MILP-II}}(\mathbf y, \mathbf s, \mathbf x, \mathbf g, \mathbf h,\mathbf{\bar h}, \mathbf z, \bm\alpha,\bm{\beta}, \bm \rho)
\label{eq:phase2_obj}
\\
\mathrm{s.t.}\quad
& (\mathbf y, \mathbf s, \mathbf x, \mathbf g, \mathbf h,\mathbf{\bar h}, \mathbf z, \bm\alpha,\bm{\beta}, \bm \rho)\in \mathcal{X}^{\mathrm{MILP-II}}.
\label{eq:phase2_original_region}
\\
& y_k\leq \hat y^L_k,
&& k\in K,
\label{eq:phase2_y_fixed}
\\
& \sum_{k\in K}\beta_{ik}+s_i=1,
&& i\in I.
\label{eq:phase2_no_rejection}
\end{align}
\end{subequations}

Constraints \eqref{eq:phase2_y_fixed} ensure that the collection points opened in Phase 2 are selected within those activated in Phase 1. Constraints \eqref{eq:phase2_no_rejection} require each customer to be assigned to exactly one service mode: either home delivery or collection-point pickup and enforces that a customer routed for home delivery cannot reject the proposed service. 
Note that these constraints are not valid inequalities for the original \ABC, in which customers are allowed to reject the proposed home-delivery option and choose a collection point instead.\\
The following constraints have been incorporated to the above formulation in order to enhance the empirical performance of the solver and to facilitate the construction of high-quality feasible solutions:  

\begin{subequations}
\begin{align}
& \sum_{\substack{r\in R^i\\r\neq \bar r^i}}\alpha_{ir}\leq s_i,
&& i\in I,
\label{eq:phase2_home_acceptance} \\
& \alpha_{ir} + \beta_{ik} \leq 1,
&& i \in I,\ r \in R^i,\ k \in K, \label{eq:phase2_slot_acceptance} \\
& && r \neq \bar r^i,\ \gamma_{ir} \geq \gamma_{i,k+|R^i|}. \notag
\end{align}
\end{subequations}


Let $(\mathbf{\hat y}, \mathbf{\hat s}, \mathbf{\hat x}, \mathbf{\hat  g}, \mathbf{\hat h},\mathbf{{\hat{\bar h}}}, \mathbf{\hat z}, \bm{\hat \alpha},\bm{\hat \beta}, \bm{\hat \rho})$  
denote an optimal solution obtained in Phase 2. This solution provides a structured feasible solution to the original problem and, therefore, an upper bound on its optimal objective value.\\ 

The resulting two-phase MIP heuristic can be summarized as follows:

\begin{enumerate}
    \item Solve the Phase 1 reduced location-routing problem \eqref{eq:phase1_obj}-\eqref{eq:phase1_prop1} to identify the promising collection-point set $K^L$.
    \item Extract the best Phase 1 solution, including the opening, routing, and customer-to-collection-point assignment decisions.
    \item Use the solution of Phase 1 to initialize the corresponding variables in Phase 2.
    \item Solve the Phase 2 restricted \ABC formulation given by \eqref{eq:phase2_obj}, \eqref{eq:phase2_home_acceptance}-\eqref{eq:phase2_slot_acceptance}.
    \item Return the solution of Phase 2 as the output of the MIP heuristic. This solution is feasible for the original problem and provides an upper bound on its optimal objective value.
\end{enumerate}

\begin{figure}[h]
\centering
\includegraphics[width=14 cm]{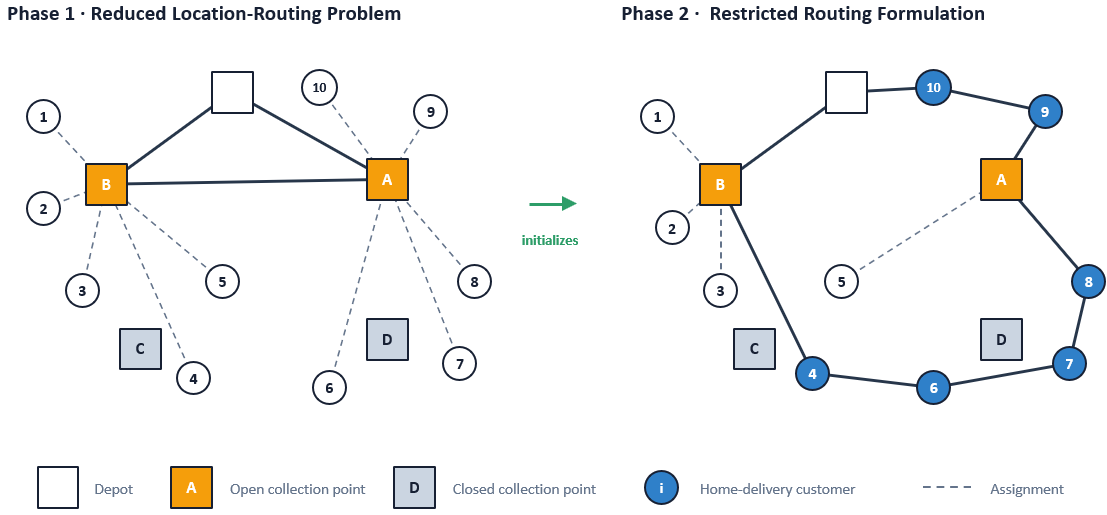}
\caption{MIP heuristic procedure.}
\label{fig:MIP_Heur}
\end{figure}

Figure \ref{fig:MIP_Heur} illustrates the two phases of the proposed MIP heuristic. In Phase 1, the reduced location-routing problem is solved to identify a promising subset of collection points. The selected collection points, shown in orange, are connected to the depot through a preliminary route, while customers are assigned to their preferred selected collection points, as represented by the dashed lines. Collection points that are not selected in this phase are shown in gray. In Phase 2, the restricted \ABC\, formulation is solved using the Phase 1 solution as an initialization. The collection-point configuration identified in Phase 1 is maintained, while the routing decisions are extended to include customers selected for home delivery, shown in blue. Since customer rejection is not permitted in this phase, each customer is either served through home delivery or assigned to one of the selected collection points. The resulting Phase 2 solution is feasible for the original problem and constitutes the final solution returned by the MIP heuristic.

\subsection{Warm-Start Strategy for the MILP Formulation}
As will be seen, the MIP heuristic described in the above section produces a feasible solution that, in general, is of good quality. Moreover, the heuristic solution can also be used as a warm-start solution to enhance the performance of the MILP formulation(s). 
Two alternative strategies can be followed: 
($i$) the solver runs for a limited (small) computing time and stops immediately after finding a feasible solution better than that of the warm start; or, ($ii$) the solver uses the heuristic solution to reduce the search space to those better than it, and runs until a provable optimal solution has been found or until a maximum (large) time limit has been reached.\\

In either case, when the heuristic is applied prior the MILP, the solution produced by the heuristic is provided to the solver as a warm-start solution, which is used as the initial incumbent solution. 
The solver then focuses the search on finding solutions better than the initial one and can either preserve or improve the initial solution. \\

\begin{figure}[h]
\centering
\includegraphics[width=14 cm]{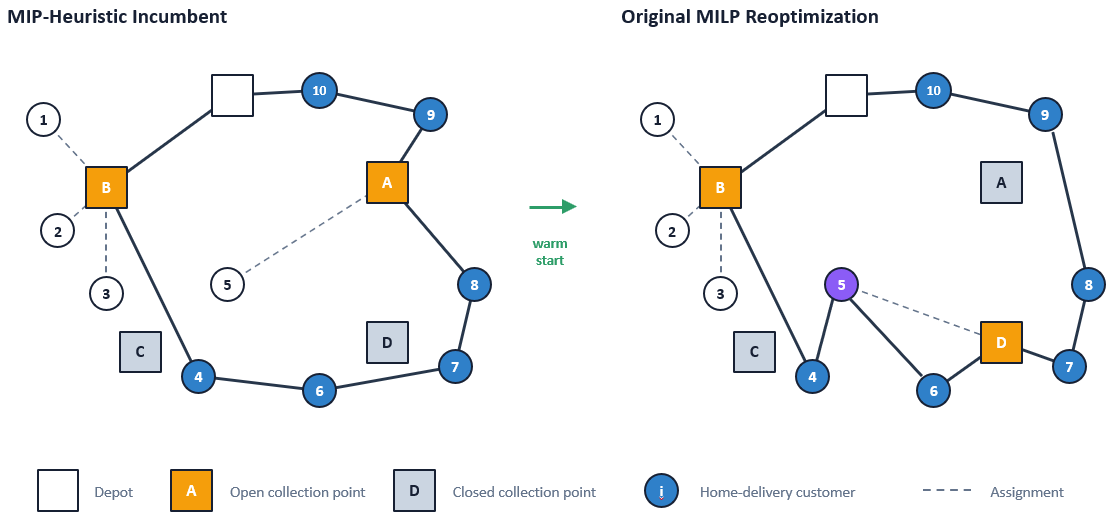}
\caption{Warm-start procedure using the MIP heuristic.}
\label{fig:MIP_Heur_warm}
\end{figure}

Figure \ref{fig:MIP_Heur_warm} illustrates how the solution obtained after reoptimization may differ from the one generated by the MIP heuristic. Several changes can be observed with respect to the Phase 2 solution shown in Figure 2. In particular, the set of activated collection points changes, leading to a different routing and assignment structure. Moreover, customer 5 is included in the route and is therefore offered home delivery; however, the customer rejects the proposed service time and ultimately chooses the pickup mode at an activated collection point. This behavior, highlighted in purple, illustrates an important feature that cannot occur in the restricted Phase 2 formulation: being included in the delivery route does not necessarily imply that the customer accepts home delivery. Overall, the figure shows how reoptimization can simultaneously revise collection- point activation, routing, and customer-assignment decisions while recovering the endogenous customer responses of the original problem.

\section{Experimental Results}
\label{sect:results}

In this section, we present the results of the computational experiments conducted for the \ABC. First, we compare the computational performance of the SLR, MILP-I, and MILP-II formulations developed for the problem under study. Based on this comparison, we then use MILP-II to assess the computational tractability of the \ABC\ across instances of increasing size and to identify the limitations of the exact solution approach. Finally, we evaluate the effectiveness of the proposed MIP heuristic, with particular emphasis on large-size instances for which the exact formulation becomes computationally challenging. All experiments were conducted on a desktop computer equipped with an Intel® Core™ i7-11370H @ 3.30 GHz processor and 32 GB of RAM. The mathematical formulations and the MIP heuristic were implemented in Python 3.10 and solved using Gurobi 13.0. Before presenting the computational results, we first describe the benchmark instances used in the experiments.

\subsection{Benchmark Instances}
Since the \ABC\ is a new problem, to the best of our knowledge, no benchmark instances are currently available in the literature for the variant introduced in this study. Therefore, we generate a new benchmark set comprising 35 instances specifically designed for the \ABC. The generation scheme follows the principles of synthetic benchmark design commonly adopted in vehicle-routing studies, in which controlled pseudo-random instances are generated to assess computational performance across different problem sizes and structural characteristics.\\

Each instance considers a single depot, a set of customers, and a set of candidate collection points. The instances are classified according to the number of customers as follows: small-size instances with $|I| \in \{20,30\}$, medium-size instances with $|I| \in \{40,50,60\}$, and large-size instances with $|I| \in \{80,100\}$. Five instances were generated for each value of $|I|$, using five different random seeds. The number of candidate collection points is determined by a density parameter $\rho$, such that $|K|=\operatorname{round}(\rho |I|)$. In all instances, $\rho=0.20$, resulting in $|K| \in \{4,6,8,10,12,16,20\}$ for $|I| \in \{20,30,40,50,60,80,100\}$, respectively.\\

Customer, collection-point, and depot coordinates are generated on a two-dimensional $3000 \times 3000$ grid. Following the methodology of \cite{uchoa2017new}, all nodes are randomly located within the service region. Travel times are computed based on Manhattan distances, which are suitable for representing movement in grid-like urban networks. Specifically, for each pair of nodes $i$ and $j$, the travel time is given by $t_{ij}=d_{ij}/\upsilon$, where $d_{ij}$ denotes the Manhattan distance between the two nodes and $\upsilon$ is the travel-speed parameter, measured in meters per minute. We assume an average operational speed of 15 km/h, corresponding to $\upsilon=250$ m/min. This value represents an effective urban delivery speed that accounts for the operational conditions typically encountered in last-mile distribution, including intersections, congestion, parking, and stop-and-go movements.\\

Routing costs are computed as $c_{ij}=\kappa\, d_{ij}/1000$ where $\kappa$ denotes the routing-cost parameter per kilometer. We set $\kappa=25$ for all instances. Rather than representing a currency-calibrated transportation cost, this value is used as a normalized cost coefficient to balance routing costs with collection-point activation, pickup, and temporal-deviation costs. Each candidate collection point has a fixed opening cost $f_k$, randomly selected from the set $\{50,55,60,65\}$. These values introduce moderate cost variation among candidate collection points while maintaining opening costs within a comparable range, thereby preventing collection point-opening costs from dominating the remaining components of the objective function.\\

The compensation cost associated with collection point $k$ being the closest activated collection point to customer $i$, who does not receive home delivery, is computed based on the Manhattan distance $d$ between them. Specifically, $\theta^F_{ik}=\varphi\, d_{ik}/1000$, where $\varphi$ denotes the pickup compensation cost per kilometer. We set $\varphi=25$, using the same scale as the routing-cost coefficient. This choice provides a comparable valuation of vehicle travel and customer pickup inconvenience, facilitating the trade-off between home delivery and a potentially preferable collection-point assignment. In addition, each customer has a preference ranking over the candidate collection points. These rankings are generated using a randomized distance-based measure, such that preferences are influenced by distance without being determined exclusively by geographic proximity.\\
 
To preserve numerical stability and avoid artifacts arising from limited numerical precision, $t_{ij}$, $d_{ij}$, $d_{ik}$, and $c_{ij}$ are rounded to two decimal places. Preliminary computational experiments indicated that this level of precision helps prevent degenerate cases, such as zero-time or zero-cost arcs, which could weaken the sequencing logic of the model. It is worth noting, however, that subtour elimination is ensured by the routing constraints of the formulation and does not rely on the adopted rounding scheme.\\

All instances use a planning horizon of $T_{\max}=480$ minutes, corresponding to an eight-hour operating period. Each customer is assigned one or two feasible home-delivery time intervals, with durations randomly selected from $\{0.5, 1\}$ hours. The preference ranks are assigned independently of their chronological order. The most preferred interval determines the customer's ideal delivery time, which is located either at the beginning or at the end of that interval. Deviations from the ideal delivery time are penalized at $\theta_i^T=0.25$ cost units per minute (15 per hour), allowing the model to balance time preferences against routing, pickup, and facility-opening costs. Table~\ref{tab:benchmark_instances} summarizes the main dimensions of the benchmark, including the number of customers $|I|$, candidate collection points $|K|$, and feasible home-delivery intervals $|R^i|$ for each customer. Five instances are generated for each $(|I|,|K|)$ configuration, yielding 35 instances in total.

\begin{table}[htbp]
\centering
\caption{Size of the benchmark instances.}
\label{tab:benchmark_instances}
\small
\setlength{\tabcolsep}{10pt}
\renewcommand{\arraystretch}{1.10}
\begin{tabular}{cccc}
\toprule
$\boldsymbol{|I|}$
& $\boldsymbol{|K|}$
& $\boldsymbol{|R^i|}$
& \textbf{\# of instances} \\
\midrule
20  & 4  & 1 or 2 & 5 \\
30  & 6  & 1 or 2 & 5 \\
40  & 8  & 1 or 2 & 5 \\
50  & 10 & 1 or 2 & 5 \\
60  & 12 & 1 or 2 & 5 \\
80  & 16 & 1 or 2 & 5 \\
100 & 20 & 1 or 2 & 5 \\
\midrule
\multicolumn{3}{r}{\textbf{Total}} & \textbf{35} \\
\bottomrule
\end{tabular}%
\end{table}

\subsection{Computational Comparison of the Formulations}

An additional set of 30 instances was generated exclusively to compare the computational performance of the SLR, MILP-I, and MILP-II formulations. For this experiment, instances with $|I|\in\{20,25\}$ customers and $|K|\in\{4,5\}$ candidate collection points were considered, while all remaining parameters were set as described in the previous subsection. The comparison is intended to assess the computational benefits provided by the two proposed MILP reformulations relative to the classical SLR. A time limit of 600 seconds was imposed for each run. Table \ref{tab:summary_all_formulations} reports the resulting computational performance. Recall that SLR denotes the classical single-level reformulation based on primal-dual optimality conditions with non-linearities, MILP-I its linearized counterpart obtained with the linearizations described in Section  \ref{sec:MILP-I}, 
and MILP-II the alternative reformulation using primal variables only, derived from the properties of the sets $K^z_{\prec}(i,r)$ and $K^y_{\succ}(i,k)$ (see Proposition \ref{propo:primal}).\\

\begin{table}[htbp] 
\centering 
\caption{Comparison of solution values and computational times for the three formulations.} 
\label{tab:summary_all_formulations} 
\begin{tabular}{llrcrr} 
\toprule 
\textbf{Size} & \textbf{Formulation} 
& \multicolumn{2}{c}{\textbf{Solution value}} 
& \multicolumn{2}{c}{\textbf{CPU time (s)}} \\ 
\cmidrule(lr){3-4} \cmidrule(lr){5-6} 
 & & \textbf{Avg.} & \textbf{Std. dev.} 
 & \textbf{Avg.} & \textbf{Std. dev.} \\ 
\midrule 
20 & SLR     & 628.06 & 27.81 & 15.75  & 14.48\phantom{ab} \\ 
20 & MILP-I  & 628.06 & 27.81 & 12.64  & 8.77\phantom{ab} \\ 
20 & MILP-II & 628.06 & 27.81 & 10.27  & 6.77\phantom{ab} \\ 
\midrule 
25 & SLR     & 721.83 & 31.15 & 231.39 & 158.18\phantom{ab} \\ 
25 & MILP-I  & 721.83 & 31.15 & 127.53 & 201.01\phantom{ab} \\ 
25 & MILP-II & 721.83 & 31.15 &  81.80 & 87.89\phantom{ab} \\ 
\bottomrule 
\end{tabular}%
\end{table}

Table \ref{tab:summary_all_formulations} shows that all three formulations yield the same objective values, confirming the validity of the proposed MILP reformulations. However, substantial differences are observed in their computational performance. Both MILP-I and MILP-II reformulations outperform the SLR in terms of solution time, with MILP-II consistently exhibiting the best computational performance. These results demonstrate the computational advantage of exploiting the properties of the sets $K^z_{\prec}(i,r)$ and $K^y_{\succ}(i,k)$ in the alternative reformulation. Based on these findings, MILP-II is adopted for the remainder of the computational experiments.

\subsection{Computational Performance of MILP-II}

In this section, we evaluate the computational performance of the MILP-II reformulation for the \ABC\, across all benchmark instances. A time limit of 3600 seconds was imposed for each instance. Table~\ref{tab:summary_mip_results} summarizes the results, where \textbf{$|I|$} and \textbf{$|K|$} denote the number of customers and candidate collection points, respectively; \textbf{Avg. OF} reports the average objective-function value; \textbf{Avg. GAP (\%)} gives the average percent optimality gap; and \textbf{Avg. time (s)} reports the average computing time. The last two columns, \textbf{Feasible} and \textbf{Optimal}, report the number of instances for which a feasible solution was found and the number solved to optimality within the solver's prescribed tolerances, respectively.

\begin{table}[htbp]
\centering
\caption{Summary of MILP-II computational results.}
\label{tab:summary_mip_results}
\begin{tabular}{rrrrrrr}
\toprule
& & & & &
\multicolumn{2}{c}{\textbf{Instances}} \\
\cmidrule(lr){6-7}
\textbf{$|I|$} &
\textbf{$|K|$} &
\textbf{Avg. OF} &
\textbf{Avg. GAP (\%)} &
\textbf{Avg. time (s)} &
\textbf{Feasible} &
\textbf{Optimal} \\
\midrule
20  & 4  &  693.01 &  0.00 &    4.22 & 5/5 & 5/5 \\
30  & 6  &  851.80 &  0.00 &   59.25 & 5/5 & 5/5 \\
40  & 8  & 1070.25 &  2.59 & 2009.77 & 5/5 & 3/5 \\
50  & 10 & 1244.42 &  7.39 & 3246.33 & 5/5 & 1/5 \\
60  & 12 & 1368.85 & 10.26 & 3600.00 & 5/5 & 0/5 \\
80  & 16 & 1797.88 & 22.96 & 3600.00 & 5/5 & 0/5 \\
100 & 20 & 2718.10 & 44.12 & 3600.00 & 5/5 & 0/5 \\
\bottomrule
\end{tabular}
\end{table}

The results show that MILP-II performs effectively on the smallest instances. All instances with 20 and 30 customers are solved to optimality, with average computing times of 4.22 and 59.25 seconds, respectively. Computational effort increases substantially for 40 customers, where three of the five instances are reported as optimal, the average runtime reaches 2009.77 seconds, and the average gap rises to 2.59\%. For 50 customers, only one instance is reported as optimal, with an average runtime of 3246.33 seconds and an average gap of 7.39\%. All instances with 60 or more customers reach the time limit, with average gaps increasing from 10.26\% for 60 customers to 44.12\% for 100 customers. Overall, MILP-II obtains feasible solutions for all 35 instances and reports 14 as optimal. Although MILP-II performs well on small instances, its ability to obtain high-quality solutions within the time limit deteriorates substantially for medium and large instances, motivating the use of the proposed MIP heuristic.

\subsection{Performance of the MIP Heuristic and Effect of the Warm-Start Solution on MILP-II}

The proposed MIP heuristic is used to generate high-quality initial solutions for large \ABC\ instances, for which MILP-II becomes computationally challenging, as shown in the previous subsection. The experiment compares the solution produced by the heuristic with those obtained with MILP-II under two different strategies. The first strategy solves MILP-II directly without an initial solution, whereas in the second one, the solution generated by the two-phase MIP heuristic is used as a warm start for MILP-II. The direct approach has a time limit of 3600 seconds, while the warm-start approach allows up to 600 seconds for the two-phase heuristic in total, followed by up to 3000 seconds for MILP-II. Table~\ref{tab:mip_heuristic_comparison} summarizes the results, where $|I|$ and $|K|$ denote the number of customers and candidate collection points, respectively. Under \textbf{Objective function}, the columns \textbf{MILP-II}, \textbf{MIP heuristic}, and \textbf{WS MILP-II} report the average objective-function values obtained by direct MILP-II, the standalone heuristic, and MILP-II with a warm start, respectively. The \textbf{Deviation (\%)} columns report the percentage differences relative to the average objective-function value of direct MILP-II, with negative values indicating better solutions.

\begin{table}[htbp]
\centering
\caption{Comparison of MILP-II, the MIP heuristic, and the warm-start approach.}
\label{tab:mip_heuristic_comparison}
\begin{tabular}{rrrrrrr}
\toprule
& &
\multicolumn{3}{c}{\textbf{Objective function}}
& \multicolumn{2}{c}{\textbf{Deviation (\%)}} \\
\cmidrule(lr){3-5}
\cmidrule(lr){6-7}
\textbf{$|I|$}
& \textbf{$|K|$}
& \textbf{MILP-II}
& \textbf{MIP}
& \textbf{WS}
& \textbf{MIP }
& \textbf{WS} \\
& & & \textbf{heuristic}
& \textbf{MILP-II}
& \textbf{heuristic}
& \textbf{MILP-II} \\
\midrule
20  & 4  &  693.01 &  701.50 &  693.01 &   1.23 &   0.00 \\
30  & 6  &  851.80 &  866.65 &  851.80 &   1.74 &   0.00 \\
40  & 8  & 1070.25 & 1063.81 & 1063.81 &  -0.60 &  -0.60 \\
50  & 10 & 1244.42 & 1239.81 & 1239.81 &  -0.37 &  -0.37 \\
60  & 12 & 1368.85 & 1354.02 & 1354.02 &  -1.08 &  -1.08 \\
80  & 16 & 1797.88 & 1684.73 & 1681.85 &  -6.29 &  -6.45 \\
100 & 20 & 2718.10 & 1912.57 & 1912.57 & -29.64 & -29.64 \\
\bottomrule
\end{tabular}
\end{table}

Table~\ref{tab:mip_heuristic_times} reports the average computing times in seconds. The \textbf{MILP-II} column corresponds to solving MILP-II directly, while \textbf{MIP heuristic} reports the total time required by Phases 1 and 2. The \textbf{WS MILP-II} column reports the additional computing time spent solving MILP-II initialized with the heuristic solution. Finally, \textbf{Total} reports the combined computing time of the MIP heuristic and the subsequent MILP-II optimization.\\

\begin{table}[htbp]
\centering
\caption{Average computing times in seconds.}
\label{tab:mip_heuristic_times}
\begin{tabular}{rrrrr}
\toprule
\textbf{$|I|$} &
\textbf{MILP-II} &
\textbf{MIP heuristic} &
\textbf{WS MILP-II} &
\textbf{Total} \\
\midrule
20  &    4.22 &   1.02 &    2.85 &    3.87 \\
30  &   59.25 &   7.54 &   52.84 &   60.38 \\
40  & 2009.77 &  65.47 & 1966.64 & 2032.11 \\
50  & 3246.33 & 378.71 & 2651.48 & 3030.19 \\
60  & 3600.00 & 508.93 & 3000.00 & 3508.99 \\
80  & 3600.00 & 600.00 & 3000.00 & 3600.00 \\
100 & 3600.00 & 600.00 & 3000.00 & 3600.00 \\
\bottomrule
\end{tabular}
\end{table}

For instances with 20 and 30 customers, the MIP heuristic produces average objective values 1.23\% and 1.74\% higher than the direct MILP-II, respectively, while the warm-start approach recovers the same average values to two decimal places. From 40 customers onward, the heuristic alone outperforms direct MILP-II on average, with improvements ranging from 0.37\% to 29.64\%. For 80 and 100 customers, these improvements are obtained after approximately 600 seconds in Phase 2, in addition to the first heuristic phase, compared with 3600 seconds for direct MILP-II. The subsequent warm-start optimization provides limited additional improvements for instances with 40 or more customers. The largest additional reduction occurs for 80 customers, where the average objective value decreases from 1684.73 to 1681.85, increasing the improvement over direct MILP-II from 6.29\% to 6.45\%. For 100 customers, the heuristic and warm-start approaches both yield an average objective value of 1912.57, which is 29.64\% lower than that obtained by direct MILP-II. Thus, for larger instances, the main benefit comes from the heuristic itself, which generates substantially better solutions early in the solution process, whereas the subsequent MILP-II optimization provides comparatively modest additional improvements.\\

\begin{figure}[htbp!]
    \centering
    \includegraphics[width=0.80\textwidth]{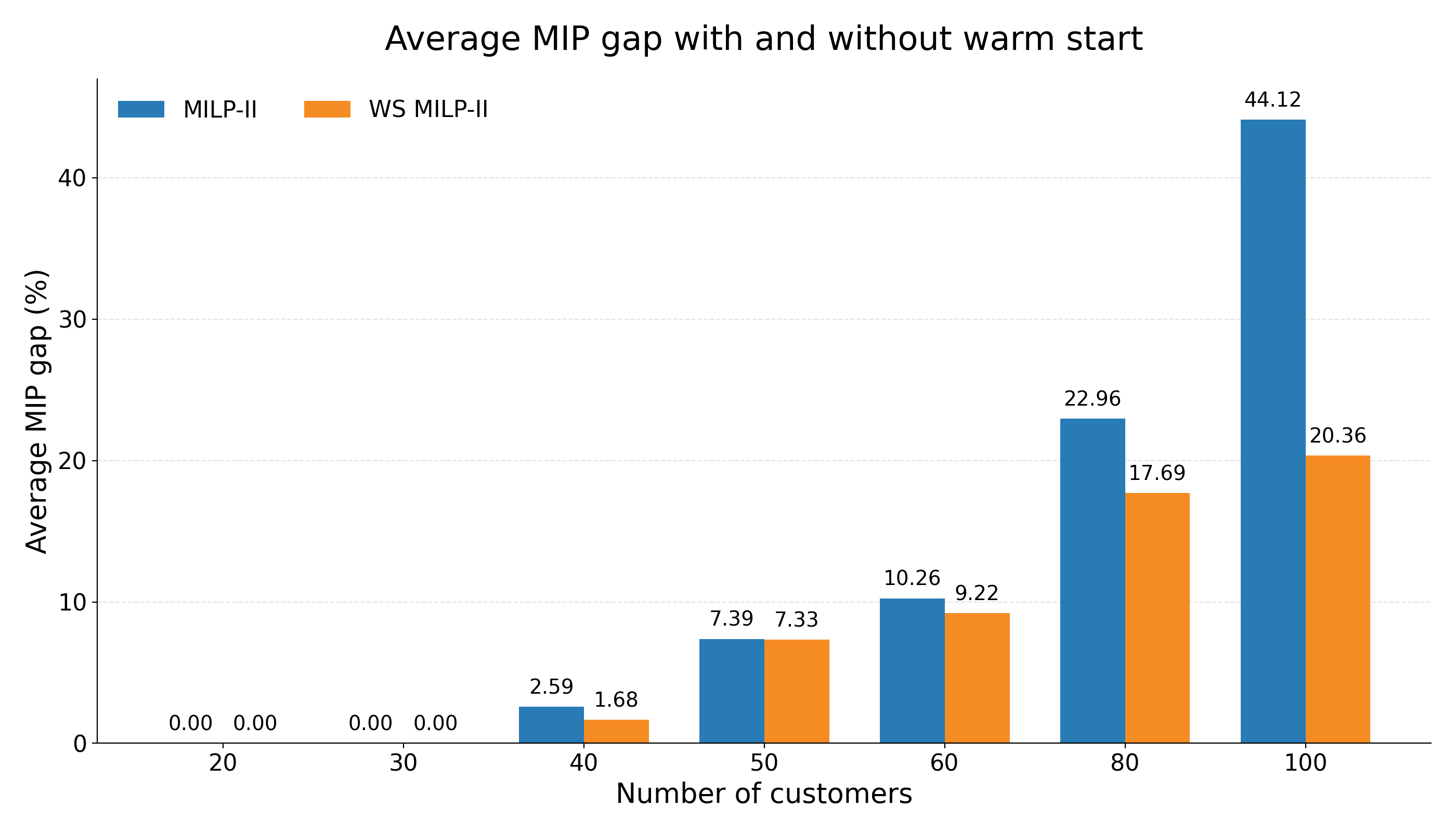}
    \caption{Average MIP gaps for pure MILP-II and the warm-start approach.}
    \label{fig:mip_gap}
\end{figure}

Figure~\ref{fig:mip_gap} compares the average percent optimality gaps obtained by MILP-II with and without the proposed warm start. The warm-start strategy reduces the average gap for instances with 40 or more customers. For 60 customers, the gap decreases from 10.26\% to 9.22\%, while for 100 customers it decreases from 44.12\% to 20.36\%, corresponding to a reduction of 23.76 percentage points. Although the gaps remain substantial for the largest instances, these reductions indicate that the warm start enables MILP-II to make greater progress within the available computational budget.

\subsection{Managerial Insights}\label{sect:managerial}
\subsubsection{Value of Customer Preferences}\label{sect:ValoCust}

\begin{figure}[h]
    \centering
    \includegraphics[width=0.99\textwidth]{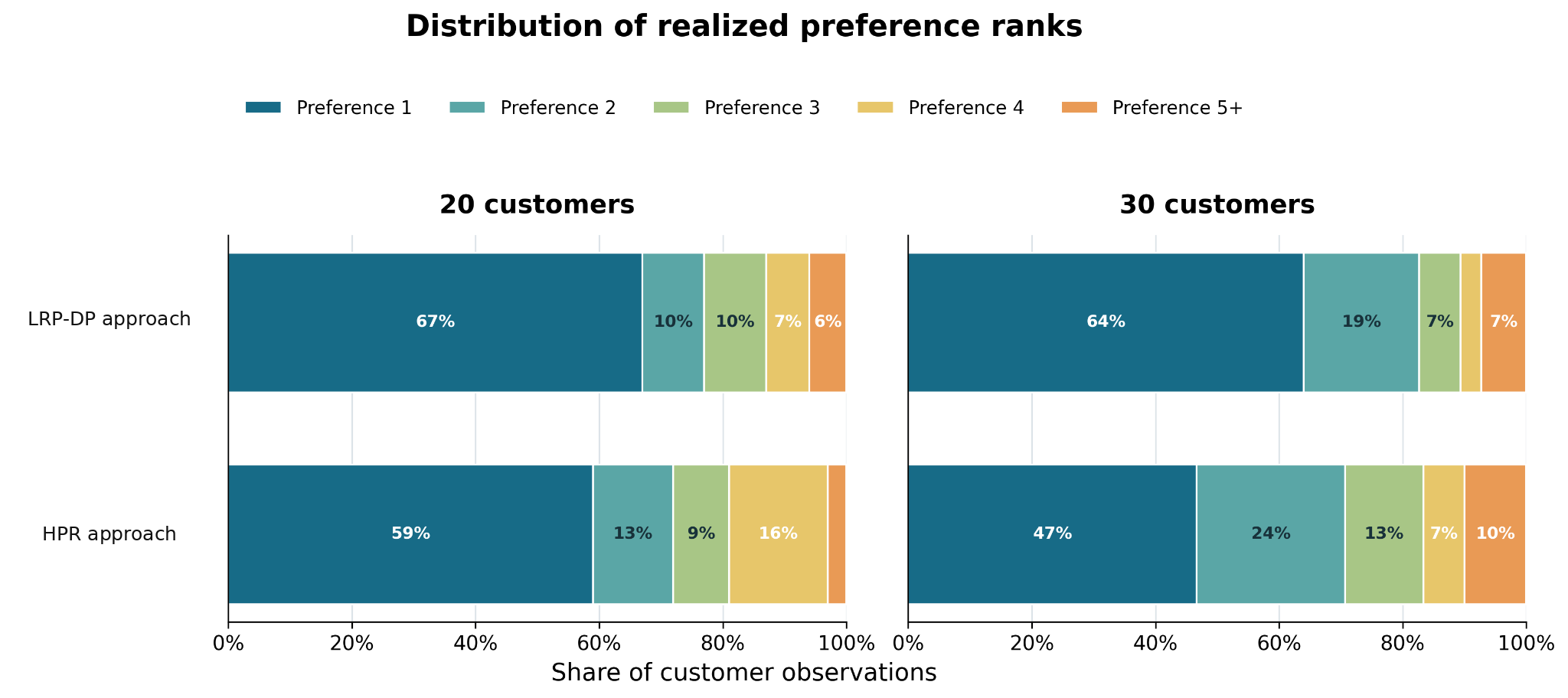}
    \caption{Preference rank distribution for the proposed bilevel framework and the leader perspective model.}
    \label{fig:prefernces}
\end{figure}

To quantify the benefit of explicitly incorporating customer preferences into the planning process, we compare the solution obtained with the proposed \ABC\ against a benchmark model that ignores customer preferences. The latter represents a setting in which the company assigns customers to collection points without allowing them to choose among the available alternatives. This approach corresponds to the high-point relaxation (HPR) commonly used in bilevel optimization \citep{bard1998practical}. In summary, both models address the same location-routing problem; however, the benchmark model omits the preference-related constraints. Under the HPR, customers are assigned to collection points by the company, so the assignment and compensation terms coincide, i.e., $\sum_{i\in I}\sum_{k\in K}\theta_{ik}^F\beta_{ik}=\sum_{i\in I}\sum_{k\in K}\theta_{ik}^F\rho_{ik}$. Note that under the HPR, customers accept the company's proposed service option, whether home delivery or pickup at a collection point. Figure~\ref{fig:prefernces} visualizes the resulting distribution of customer preference ranks. For instances with 20 customers, the percentage receiving their first-ranked option increases from 59\% under the HPR modeling approach to 67\% under the \ABC\ modeling approach; for instances with 30 customers, it increases from 47\% to 64\%. The average preference rank decreases from 1.91 to 1.77  and from 2.11 to 1.74, respectively (i.e., $-7.33\%$ and $-17.53\%$, respectively), with lower values indicating more preferred service alternatives.  These results show that incorporating customer preferences impacts location-routing decisions and improves preference satisfaction. The average optimal objective value increases from 657.81 under HPR to 669.84 under \ABC\ for 20 customers (+1.83\%), and from 827.58 to 863.66 for 30 customers (+4.36\%). These increases quantify the additional cost incurred by the company to achieve the higher level of customer preference satisfaction observed under the \ABC. 
In particular, for the considered instances, the \ABC\ produces a percent improvement in the average preference rank of $7.33\%$ at a cost increase of $1.83\%$ for $n=20$, and a percent improvement of $17.53\%$ at a cost increase of $4.36\%$ for $n=30$.\\

Overall, the managerial insights presented in this section illustrate the benefits of explicitly incorporating customer preferences into the model. Compared with the HPR modeling approach, the \ABC\ improves preference satisfaction by increasing the percentage of customers receiving their first-ranked option and reducing the average preference rank for both instance sizes. For the 30-customer instances, it also yields a more home-delivery-oriented configuration, with fewer pickup assignments and fewer activated collection points on average. These findings indicate that explicitly accounting for customer choice enables the company to design a service configuration that better reflects customer preferences while balancing routing, collection-point activation, and compensation costs.

\subsubsection{Solution Analysis}\label{sect:sol-analysis}
\begin{figure}[h]
    \centering
    \includegraphics[width=0.85\textwidth]{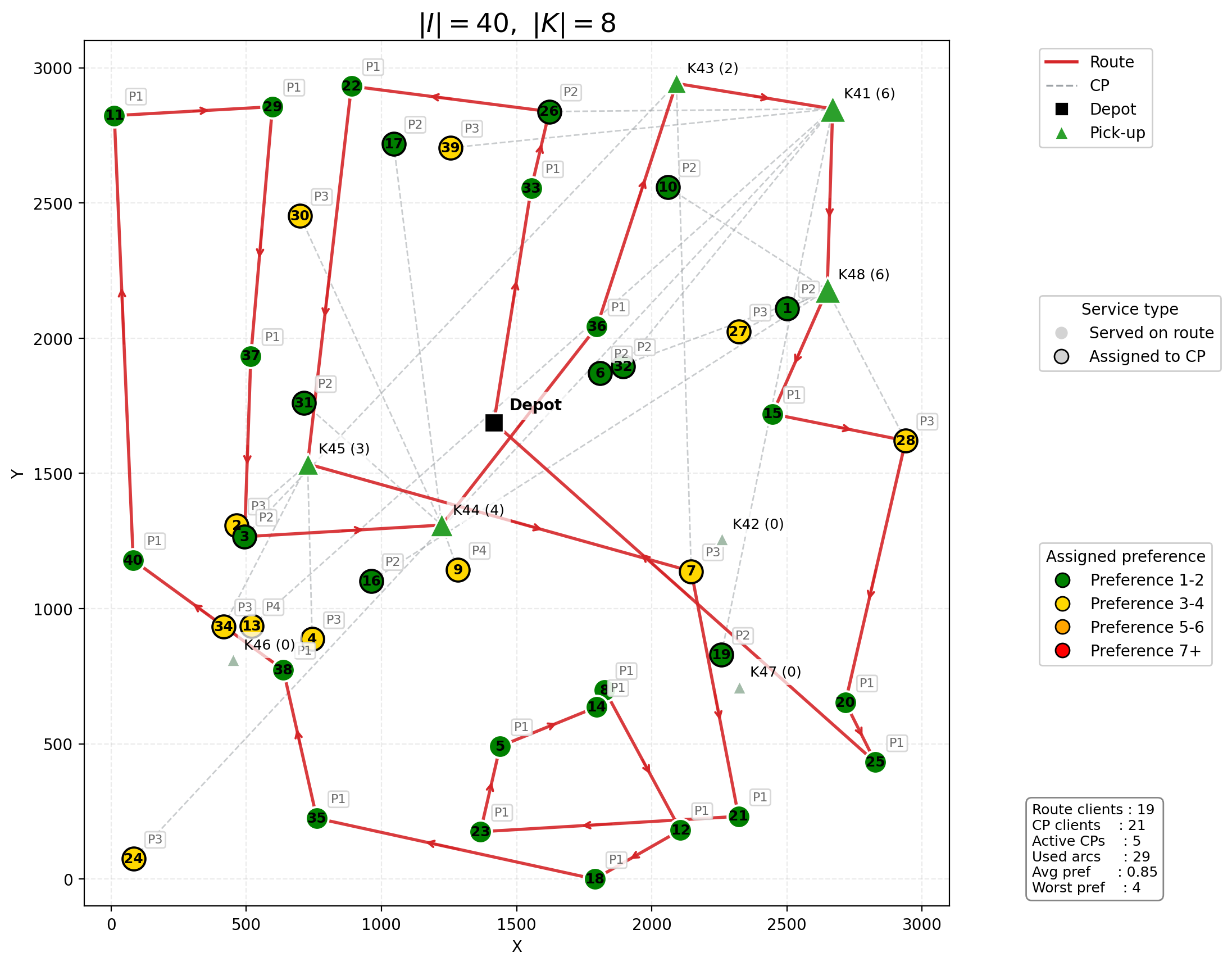}
    \caption{Feasible solution for the instance with $|I|=40$ and $|K|=8$.}
    \label{fig:seed74_feasible}
\end{figure}

\begin{figure}[h]
    \centering
    \includegraphics[width=0.85\textwidth]{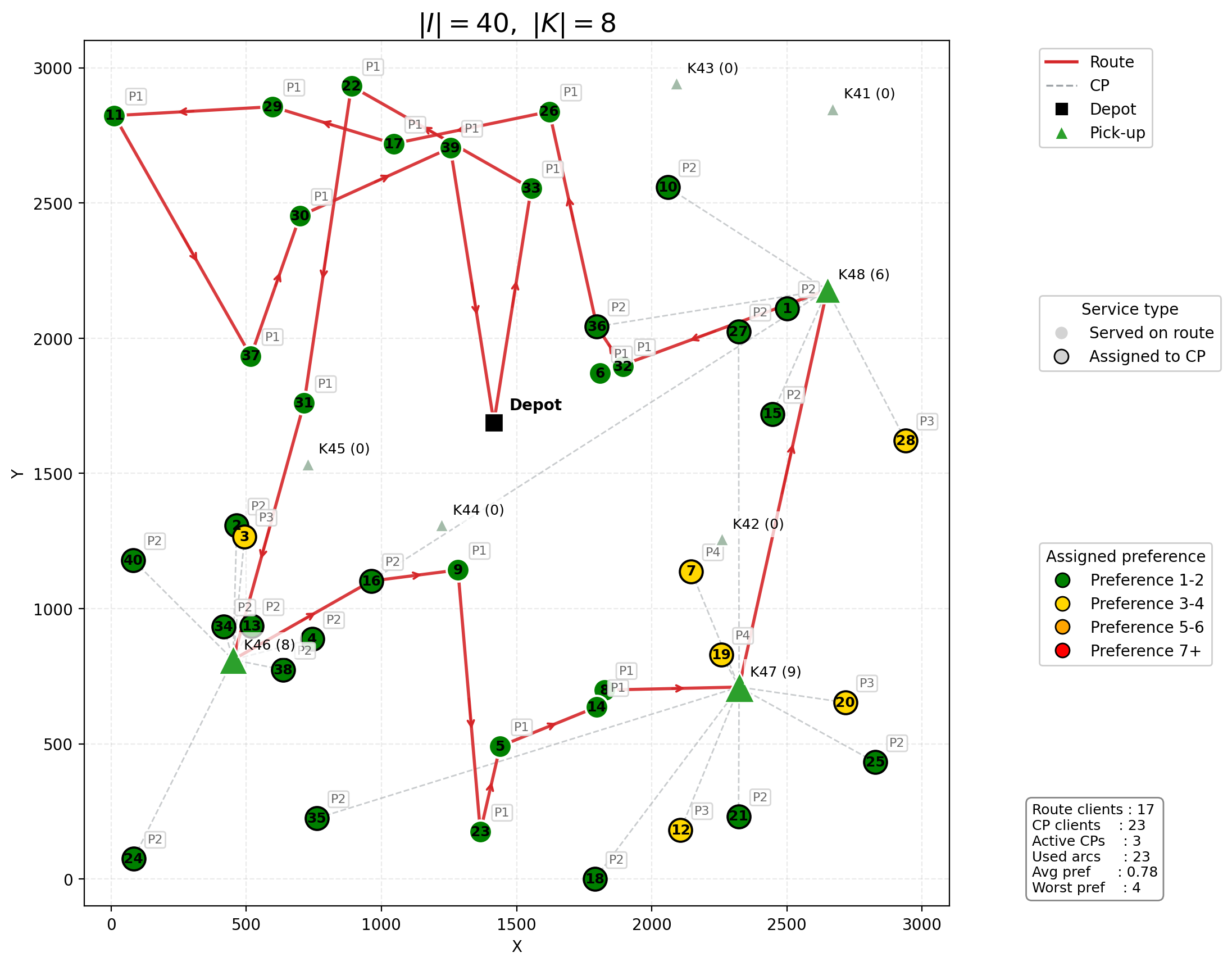}
    \caption{Optimal solution for the instance with $|I|=40$ and $|K|=8$.}
    \label{fig:seed74_optimal}
\end{figure}

In this section, we analyze an incumbent solution obtained during Gurobi's branch-and-cut search and compare it with an optimal solution to the \ABC. Both solutions correspond to the same instance, with $|I|=40$ customers and $|K|=8$ candidate collection points. The incumbent and optimal solutions have objective values of 1419.16 and 1029.07, respectively. Figures~\ref{fig:seed74_feasible} and~\ref{fig:seed74_optimal} illustrate the corresponding configurations. This comparison provides further insight into how routing and location decisions interact with customer preference satisfaction. \\

Figures~\ref{fig:seed74_feasible} and~\ref{fig:seed74_optimal} illustrate the flexibility available to customers when deciding whether to accept the home-delivery service offered through the planned route or collect their parcels at a collection point. In the feasible solution, customers 3, 7, 26, and 28 reject the proposed home delivery and choose an activated collection point. Notably, these customers remain in the planned route, illustrating that route membership does not necessarily imply acceptance of home delivery. A similar behavior is observed in the optimal solution, where customers 16 and 36 reject home delivery and choose a collection point. Still, the optimal route visits these unserved customers. This happens because short-cutting from their predecessors to their successors in the route would change the preference of the time interval when the successors receive home delivery and thus increase their time-related compensation. These cases show that planned home delivery does not necessarily coincide with the customer's final service choice, highlighting the trade-off between operational decisions and customer preferences. \\

To gain further insight into the solutions, Table~\ref{tab:solution_comparison} compares the main performance and service indicators of the optimal and feasible solutions of this instance. The first part of the table reports the objective value, i.e., the total cost, together with a breakdown of its main components, including collection-point activation, routing, time-related compensation, and pickup compensation costs. For these components, a negative difference indicates a higher cost in the feasible solution. The feasible solution has an objective value 37.91\% higher than the optimal one. The second part reports the number of activated collection points; the number and percentage of customers served through home delivery and collection points; the average preference rank, both overall and by service mode; and the worst rank in the merged preference list.\\

The feasible solution activates two additional collection points and serves two fewer customers through home delivery. The five activated collection points in the feasible solution serve six, two, four, three, and six customers, respectively. In contrast, the optimal solution serves more customers at home and operates three collection points, serving six, eight, and nine customers, respectively. The cost breakdown shows that the feasible solution incurs 208.86 additional units in routing cost, 109.00 additional units in collection-point activation cost, and 404.24 additional units in pickup compensation cost, while its time-related compensation cost is 0.82 units lower. Overall, these additional costs lead to a substantially higher objective value for the feasible solution. Moreover, the feasible solution uses 29 route arcs, compared with 23 in the optimal solution, which is consistent with its higher routing cost.\\

\begin{table}[h]
\centering
\caption{Comparison of the feasible and optimal solutions for the instance with $|I|=40$ and $|K|=8$.}
\label{tab:solution_comparison}
\small
\setlength{\tabcolsep}{5pt}
\renewcommand{\arraystretch}{1.10}
\resizebox{\textwidth}{!}{%
\begin{tabular}{lrrr}
\toprule
\textbf{Indicator}
& \textbf{Optimal}
& \textbf{Feasible}
& \textbf{Difference} \\
\midrule
Total cost
& 1029.07
& 1419.16
& $-37.91\%$\\
Opening cost
& 170.00
& 280.00
& $-110.00$\\
Routing cost
& 464.07
& 672.93
& $-208.86$\\
Time-related compensation cost
& 41.28
& 40.46
& $+0.82$\\
Pickup compensation cost
& 452.37
& 856.61
& $-404.24$\\
\midrule
Activated collection points
& 3
& 5
& $-2$ points\\
Customers served at home
& 17 (42.5\%)
& 19 (47.5\%)
& $-2$ customers\\
Customers assigned to collection points
& 23 (57.5\%)
& 21 (52.5\%)
& $+2$ customers\\
Average preference index (all customers)
& 0.78
& 0.85
& $-0.08$\\
Average preference index (home delivery)
& 0.00
& 0.00
& $0.00$\\
Average preference index (pickup)
& 1.35
& 1.62
& $-0.27$\\
Worst preference rank
& 3
& 3
& 0 (no change) \\
Customers that reject home delivery
& 2
& 4
& $-2$ customers\\
\bottomrule
\end{tabular}%
}
\end{table}

The optimal solution yields a lower overall average preference index, decreasing from 0.85 in the feasible solution to 0.78. The average preference index for home-delivery customers is 0.00 in both solutions, while the corresponding values for pickup customers are approximately 1.35 and 1.62, respectively. Both solutions have a worst preference rank of 3, indicating that customer-level trade-offs persist even at optimality. Moreover, the optimal solution reduces the number of customers rejecting the proposed home-delivery service from four to two. 

\section{Conclusions and Future Research}\label{sect:conclu}

This paper introduced the \ABC, a last-mile location-routing problem that explicitly incorporates service flexibility by allowing customers to choose between home delivery and activated collection points according to their ranked preferences. The company determines which collection points to activate, which customers to offer home delivery, and the corresponding delivery route, while customers subsequently decide whether to accept the proposed home-delivery service or collect their parcels from an activated collection point. This hierarchical interaction was modeled through a bilevel optimization framework. The model also incorporates compensation mechanisms associated with service inconvenience. In particular, compensation for home delivery is based on the deviation of actual service times relative to the customers ideal service times, whereas compensation for collection-point service is based on the closest activated collection point rather than on the collection point ultimately selected by the customer, thereby separating the service conditions provided by the company from additional travel resulting from the customer's individual preferences.\\

From a methodological perspective, we analyzed the structure of the followers' problems and exploited it to derive equivalent single-level reformulations of the bilevel model. Two MILP reformulations were proposed, with MILP-II, which uses only primal variables, providing the best computational performance among the formulations considered. Nevertheless, the computational experiments also showed that the problem becomes increasingly challenging as the number of customers and candidate collection points grows. To address this limitation, we developed a two-phase MIP heuristic that constructs feasible solutions and uses them to warm-start MILP-II. The resulting strategy becomes particularly effective for the larger instances, improving both the final objective values and the optimality gaps within the same computational budget.\\

The computational analysis also highlights the operational relevance of explicitly considering customer choice. In particular, comparison with the high-point relaxation, in which the company controls the collection-point assignments, shows that accounting for customer preferences changes the resulting location-routing decisions and improves preference satisfaction. The bilevel solutions increase the percentage of customers receiving their first-ranked service alternative and reduce the average preference rank. The results show that the increase in the overall cost is moderate compared with the improvement in customer satisfaction. These results indicate that ignoring customer autonomy may lead to service configurations that appear attractive from the company's perspective but are less aligned with the service that customers would actually choose. Thus, incorporating customer choice provides a more realistic representation of the interaction between operational planning and service preferences in last-mile delivery.\\

Several extensions provide opportunities for future research. First, the computational difficulty observed for larger instances motivates the development of more scalable exact and heuristic algorithms. A related direction is the incorporation of machine learning techniques to enhance the performance of the proposed optimization methods. Second, the models assume deterministic customer preference rankings. Considering uncertain, partially observed, or dynamically evolving preferences would provide a natural extension and could lead to stochastic or robust variants of the \ABC. Finally, the problem could be extended by incorporating additional features commonly considered in location-routing problems, such as capacitated collection points, heterogeneous delivery resources, and other operational constraints.

\section*{Disclosure Statement}
No potential conflict of interest was reported by the authors.

\section*{Funding}
J.F. Camacho-Vallejo acknowledges financial support by grant RED2022-134149-T funded by MICIU/AEI/10.13039/501100011033 (Thematic Network on Location Science and Related Problems). The research of E. Fernández has been partially funded by the Spanish Ministry of Science and Innovation through project PID2023-146643NB-I00 of the State Research Agency (AEI).

\bibliographystyle{elsarticle-harv}\biboptions{authoryear}
\bibliography{Bibliography}

\end{document}